\documentclass[graybox]{svmult}

\usepackage{type1cm}        
\usepackage{makeidx}         
\usepackage{graphicx}        
\usepackage{multicol}        
\usepackage[bottom]{footmisc}

\usepackage{newtxtext}       %
\usepackage{newtxmath}       

\makeindex             
\usepackage{url,doi,hyperref}

\let\theoremstyle\relax

\let\proof\relax

\usepackage{amsmath, amsthm, amsfonts}

\usepackage{hyperref}
\usepackage{epic}
\usepackage[dvips]{epsfig}
\usepackage{enumitem}
\usepackage{csquotes}

\theoremstyle{plain} 
  \newtheorem{theo}{Theorem}[section]

\theoremstyle{definition}

  \newtheorem{assumption}[theo] {Assumption}

\newcommand{\R}{{\mathbb{R}}}

\newcommand{\calA}{{\cal A}}
\newcommand{\calB}{{\cal B}}
\newcommand{\calC}{{\cal C}}
\newcommand{\calD}{{\cal D}}

\newcommand{\calX}{{\cal X}}
\newcommand{\calY}{{\cal Y}}
\newcommand{\calU}{{\cal U}}
\newcommand{\calW}{{\cal W}}
\newcommand{\utransf}{U}

\def\bF{\mathbf{F}}
\def\by{\mathbf{y}}

\def\ttheta{\tilde{\theta}}
\def\tu{\tilde{u}}
\def\tv{\tilde{v}}
\def\cL{\mathcal{L}}

\def\R{\mathbb{R}}
\def\ts{{\tilde{s}}}
\def\tts{{\hat{s}}}
\def\ttts{{\check{s}}}

\def\dupair#1{\left\langle #1 \right\rangle}

\newcommand*{\QED}{\hfill\ensuremath{\square}}
\setlist{itemsep=5pt}

\newcommand\mpar[1]{}
\def\T{\mathbb{T}}

\begin{document}

\title*{The tangential cone condition for some coefficient identification model problems in parabolic PDEs}
\titlerunning{The tangential cone condition for parabolic inverse problems}
\author{Barbara Kaltenbacher, Tram Thi Ngoc Nguyen and Otmar Scherzer}
\authorrunning{Kaltenbacher, Nguyen and Scherzer}
\institute{Barbara Kaltenbacher \at 
Department of Mathematics,
Alpen-Adria-Universit\"at Klagenfurt, 
Universit\"atsstra\ss e 65-67, 
9020 Klagenfurt, Austria, \email{barbara.kaltenbacher@aau.at}
\and Tram Thi Ngoc Nguyen \at 
Department of Mathematics,
Alpen-Adria-Universit\"at Klagenfurt, 
Universit\"atsstra\ss e 65-67, 
9020 Klagenfurt, Austria, \email{tram.nguyen@aau.at}
\and Otmar Scherzern \at 
Faculty of Mathematics, 
University of Vienna, 
Oskar-Morgenstern Platz, 1090 Vienna, Austria, 
and Johann Radon Institute for Computational and Applied Mathematics (RICAM), 
Altenbergerstra\ss e 69, 
4040 Linz, Austria, \email{otmar.scherzer@univie.ac.at}}
%
%
\maketitle

\abstract{The tangential condition was introduced in \cite{HNS95} as a
sufficient condition for convergence of the Landweber iteration for
solving ill--posed problems.
\\
In this paper we present a series of 
time dependent
benchmark inverse problems for
which we can verify this condition.
}

\section{Introduction}

We consider the problem of recovering a parameter $\theta$ in the evolution equation
\begin{align}
& \dot{u}(t) = f(t,\theta,u(t)) \quad t\in (0,T) 
\label{model-1}\\ 
& u(0)=u_0, \label{model-2}
\end{align}
where for each $t\in(0,T)$ we consider $u(t)$ as a function on a bounded {$C^{1,1}$} domain $\Omega \subset \R^d$. In (\ref{model-1}), $\dot{u}$ denotes the first order time derivative of $u$ and $f$ is a nonlinear function. 
{
We here focus on the setting of $\theta$ not being time dependent. Problems with state equations of the form $\dot{u}(t) = f(t,\theta(t),u(t))$ could be treated analogously but this would lead to different requirements on the underlying function spaces.
}
These model equations are equipped with additional data obtained from continuous observations over time
\begin{align}
& y(t)=\calC(t,u(t)), 
\label{model-3}
\end{align} 
with 
{
an observation operator
}
$\calC$, which will be assumed to be linear; in particular, in most of what follows $\calC$ is the continuous embedding $V\hookrightarrow Y$, with $V$ and $Y$ introduced below. 

\medskip

While formulating the requirements and results first of all in this general framework, we will also apply it to a number of examples as follows.

\subsubsection*{Identification of a potential}

We study the problem of identifying the space-dependent parameter $c$ from observation of the state $u$ 
in $\Omega\times(0,T)$ 
in
\begin{alignat} {3}
&\dot{u}-\Delta u  + cu=\varphi \qquad &&(t,x) \in (0,T)\times \Omega& \label{c-problem-1}\\
&u_{|\partial\Omega}=0 && t \in (0,T)&	\label{c-problem-2}\\
& u(0)=u_0 && x \in \Omega, &	\label{c-problem-3}
\end{alignat}
where $\varphi\in L^2(0,T;H^{-1}(\Omega))$ and $u_0\in L^2(\Omega)$ are known.
Here, $-\Delta$ could be replaced by any linear elliptic differential operator with smooth coefficients.\\
With this equation, known, among others, as diffusive Malthus equation \cite{Malthus}, one can model the evolution of a population $u$ with diffusion and with exponential growth as time progresses. The latter phenomenon is quantified by the growth rate $c$, which, in this particular case, depends only on the environment.\\

\subsubsection*{Identification of a diffusion coefficient}

We further consider the problem of recovering the 
{
space-dependent
}
parameter $a$ from measurements of $u$ 
in $\Omega\times(0,T)$, 
governed by the diffusion equation
\begin{alignat}{3} 
&\dot{u}-\nabla\cdot\Bigl(a\nabla u\Bigr)=\varphi \qquad &&(t,x) \in (0,T)\times \Omega& \label{a-problem-1}\\
&u_{|\partial\Omega}=0 && t \in (0,T)&	\label{a-problem-2}\\
&u(0)=u_0 && x \in \Omega, &	\label{a-problem-3}
\end{alignat}
where $\varphi\in L^2(0,T;H^{-1}(\Omega))$ and $u_0\in L^2(\Omega)$ are known.
This is, for instance, a simple model of groundwater flow, whose temporal evolution is driven by the divergence of the flux $-a\nabla u$ and the source term $\varphi$. The coefficient $a$ represents the diffusivity of the sediment and $u$ is the piezometric head \cite{Hanke}.\\
Banks and Kunisch \cite[Chapter I.2]{BanksKunisch} discussed the more general model: $\dot{u} + \nabla\cdot\Bigl(-a\nabla u + 
b u\Bigr)+cu$, describing the sediment formation in lakes and deep seas, in particular, the mixture of organisms near the sediment-water interface.\\

\subsubsection*{An inverse source problem with a quadratic first order nonlinearity}

Here we are interested in the problem of identifying the space-dependent source term $\theta$ from observation of the state $u$
in $\Omega\times(0,T)$ 
\begin{alignat}{3} 
&\dot{u}-\Delta u - |\nabla u|^2=\theta \qquad &&(t,x) \in (0,T)\times \Omega \label{log-problem-1}\\
& u_{|\partial\Omega}=0 && t \in (0,T)	\label{log-problem-2}\\
& u(0)=u_0 && x \in \Omega\,. 	\label{log-problem-3}
\end{alignat}
This sort of PDE with a quadratic nonlinearity in $\nabla u$ arises, e.g., in stochastic optimal control theory \cite[Chapter 3.8]{FlemingSoner}.\\

\subsubsection*{An inverse source problem with a cubic zero order nonlinearity}

The following nonlinear reaction-diffusion equation involves determining the space-dependent source term $\theta$ from observation of the state $u$
in $\Omega\times(0,T)$, 
in a semiliear parabolic equation 
\begin{alignat}{3} 
&\dot{u}-\Delta u + \Phi(u)=\varphi-\theta \qquad &&(t,x) \in (0,T)\times \Omega& \label{u3-problem-1}\\
& u_{|\partial\Omega}=0 && t \in (0,T)&	\label{u3-problem-2}\\
&u(0)=u_0 && x \in \Omega, &	\label{u3-problem-3}
\end{alignat}
where the possibly space and time dependent source term $\varphi\in L^2(0,T;H^{-1}(\Omega))$ and the initial data $u_0\in H_0^1(\Omega)$ are known.

Here we selectively mention some applications for PDEs with with cubic nonlinearity $\Phi(u)$:
\begin{enumerate}
\item[] $\Phi(u)=u(1-u^2)$: 
Ginzburg-Landau equations of superconductivity \cite{BronsardStoh}, Allen-Cahn equation for the phase separation process in a binary metallic alloy \cite{AllenCahn, Nepomnyashchy}, Newell-Whitehead equation for convection of fluid heated from below \cite{GildingKersner}.

\item[] $\Phi(u)=u^2(1-u)$:
Zel'dovich equation in combustion theory \cite{GildingKersner}.

\item[] $\Phi(u)=u(1-u)(u-\alpha), 0<\alpha<1$: Fisher's model for population genetics \cite{Pao}, Nagumo equation for bistable transmission lines in electric circuit theory \cite{NagumoYoshizawaArinomoto}. 
\end{enumerate}

In part of the analysis we will also consider an additional gradient nonlinearity $\Psi(\nabla u)$ in the PDE, cf. \eqref{nonlin} below.

\medskip

Coming back to the general setting \eqref{model-1}--\eqref{model-3} we will make the following assumptions, where all the considered examples fit into.
The operators defining the model and observation equations above are supposed to map between the function spaces
\begin{align}
& f: (0,T)\times \calX\times V \rightarrow W^* \qquad\qquad\qquad\qquad\qquad\\
& \calC: (0,T)\times V \rightarrow Y,
\end{align}
where $\calX, Y, W, V \subseteq Y$ are Banach spaces.
{
More precisely, $\calX$ is the parameter space, $Y$ the data space, $W^*$ the space in which the equations is supposed to hold and $V$ the state space. The latter three are first of all the spaces for the respective values at fixed time instances and will also be assigned a version for time-dependent functions, denoted by calligraphic letters. So $\calX$, $\calY$, $\calW^*$ will denote the parameter, data and equation spaces, respectively, and $\calU$ or $\tilde{\calU}$ (to distinguish between the different versions in the reduced and all-at-once setting below) the state space.
}
The initial condition $u_0\in H$, where $H$ is a Banach space as well, will in most of what follows be supposed to be independent of the coefficient $\theta$ here. 
Dependence of the initial data and also of the observation operator on $\theta$ can be relevant in some applications but leads to further technicalities, thus for clarity of exposition we shift consideration of these dependencies to future work. 

For fixed $\theta$, we assume that the Caratheodory mappings $f$ and $\calC$ as defined above induce Nemytskii operators \cite[Section 4.3]{Troltzsch} (for which we will use the same notation $f$ and $\calC$) on the function space 
\[
\mathcal{U}=L^2(0,T;V)\cap H^1(0,T;W^*) \mbox{ or } \tilde{\calU}=L^\infty(0,T;V)\cap H^1(0,T;W^*)\,,
\]
{
cf. \eqref{Uaao} and \eqref{statespaceRe} respectively,
}
in which the state $u$ will be contained, and map into the image space $\calW^*$ and observation space $\calY$, respectively, where
\begin{equation}
\calW^*=L^2(0,T;W^*), \qquad \calY=L^2(0,T;Y). \qquad\ \label{spaceYW}
\end{equation}
Moreover, $\tilde{\calU}$ or $\calU$, respectively, will be assumed to continuously embed into $C(0,T;H)$ in order to make sense out of \eqref{model-2}.

\medskip

We will consider formulation of the inverse problem on one hand in a classical way, as a nonlinear operator equation 
\begin{equation}\label{Fthetay}
F(\theta)=y
\end{equation}
with a forward operator $F$ mapping between Banach spaces $\mathcal{X}$ and $\mathcal{Y}$,
and on the other hand also, alternatively, as a system of model and observation equation
\begin{eqnarray}
\label{Athetau}
\calA(\theta,u) &=& 0;\\
\label{Cuy}
\calC(u) &=& y.  
\end{eqnarray} 
Here, 
\begin{equation}\label{AC}
\begin{aligned}
&\calA: \mathcal{X} \times \mathcal{U} \rightarrow \mathcal{W}^*\times H, && (\theta,u)\mapsto \calA(\theta,u)=(\dot{u}-f(\theta,u),u(0)-u_0)\\
&\calC: \mathcal{U} \rightarrow \mathcal{Y}
\end{aligned}
\end{equation}
are the model and observation operators, so that with the parameter-to-state map $S:\mathcal{X} \to \mathcal{U}$ defined by 
\begin{equation}\label{S}
\calA(\theta,S(\theta))=0
\end{equation}
and 
\begin{equation}\label{FCS}
F=\calC\circ S,
\end{equation}
\eqref{Fthetay} is equivalent to the all-at-once formulation \eqref{Athetau}, \eqref{Cuy}. 
Defining 
\[
\bF:\mathcal{X} \times \mathcal{U} \rightarrow \mathcal{W}^*\times H\times \mathcal{Y}
\] 
by 
\[
\bF(\theta,u)=(\calA(\theta,u),\calC(u)), 
\]
and setting $\by=(0,y)$, we can rewrite \eqref{Athetau}, \eqref{Cuy} analogously to \eqref{Fthetay}, as 
\begin{equation}\label{bFthetauy}
\bF(\theta,u)=\by\,.
\end{equation}  
All-at-once approaches have been studied for PDE constrained optimization in, e.g., \cite{KunischSachs,KupferSachs,LeibfritzSachs,LeHe16,orozco-ghattas-97b,shenoy-heinkenschloss-cliff-98,Taasan91} and more recently, for ill-posed inverse problems in, e.g., \cite{BurgerMuehlhuberIP,BurgerMuehlhuberSINUM,HaAs01,aao16,KKV14b,LeHe16}, particularly for time dependent models in \cite{Kaltenbacher,Nguyen}.
{
The fact that we are actually using different state spaces $\calU$, $\tilde{\calU}$ in these two settings is on one hand due to the requirements arising from the need for well-definedness and differentiability of the parameter-to-state map in the reduced setting. On the other hand, while these constraints are not present in the all-at-once setting and a quite general choice of the state space is possible there, whenever a Hilbert space setting is required --- e.g., for reasons of easier implementation --- this does not only apply to the parameter and data space but also to the state and equation spaces in the all-at-once setting, whereas in a reduced setting these spaces are ``hidden'' inside the forward operator.
}

\medskip

Convergence proofs of iterative regularization methods for solving \eqref{Fthetay} (and likewise \eqref{bFthetauy}) such as the Landweber iteration \cite{HNS95,KalNeuSch08} or the iteratively regularized Gauss-Newton method \cite{BakKok04,KalNeuSch08,KaltenbacherPreviatti18} require structural assumptions on the nonlinear forward operator $F$ such as the tangential cone condition \cite{Sche95}
\begin{equation}\label{tcc}
\|F(\theta)-F(\ttheta)-F'(\theta)(\theta-\ttheta)\|_{\mathcal{Y}}\leq c_{tc}\|F(\theta)-F(\ttheta)\|_{\mathcal{Y}} 
\quad\forall \theta,\ttheta \in B_\rho^{\mathcal{X}}(\theta^0)\,,
\end{equation}
{
for some sufficiently small constant $c_{tc}$.
}
Here $F'(\theta)$ does not necessaritly need to be the Fr\'{e}chet or G\^{a}teaux derivative of $F$, but it is just required to be some linear operator that is uniformly bounded in a neighborhood of the initial guess $\theta_0$, i.e., $F'(\theta)\in L(\mathcal{X}, \mathcal{Y})$ such that    
\begin{equation}\label{Fpbounded}
\|F'(\theta)\|_{L(\mathcal{X},\mathcal{Y})}\leq C_F
\quad\forall \theta,\ttheta \in B_\rho^{\mathcal{X}}(\theta^0)\,,
\end{equation}
{
for some $C_F>0$.
}

The conditions \eqref{tcc} and \eqref{Fpbounded} 
{
enforce certain local convexity conditions of the residual $\theta\mapsto\|F(\theta)-y\|^2$, cf.\cite{Kindermann:2017}. 
}
In this sense, the conditions are structurally similar to conditions used in the analysis of Tikhonov regularization, such as those in \cite{ChKu96}. The tangential cone condition eventually guarantees convergence to the solution of \eqref{Fthetay} by a gradient descent method for the residual (and also for the Tikhonov functional). Therefore it ensures that the iterates are not trapped in local minima.

\medskip

The key contribution of this chapter is therefore to establish \eqref{tcc}, \eqref{Fpbounded} in the reduced setting \eqref{Fthetay} as well as its counterpart in the all-at-once setting \eqref{bFthetauy} for the above examples (as well as somewhat more general classes of examples) of parameter identification in initial boundary value problems for parabolic PDEs represented by \eqref{model-1}, \eqref{model-2}. In the reduced setting this also involves the proof of well-definedness and differentiability of the parameter-to-state map $S$, whereas in the all-at-once setting this is not needed, thus leaving more freedom in the choice of function spaces. Correspondingly, the examples classes considered in Section \ref{sec:aao} will be more general than those in Section \ref{Reduced}.

Some non-trivial static benchmark problems where the tangential condition has been verified can be found e.g., in \cite{DunHoh14,HubSheNeuSch18,LecRie08}.

\medskip

We mention in passing that in view of existing convergence analysis for such iterative regularization methods for \eqref{Fthetay} or \eqref{bFthetauy} in rather general Banach spaces we will formulate our results in general Lebesgue and Sobolev spaces. Still, we particularly strive for a full Hilbert space setting as preimage and image spaces $X$ and $Y$, since derivation and implementation of adjoints is much easier then, and also the use of general Banach spaces often introduces additional nonlinearity or nonsmoothness.  
Moreover we point out that while in the reduced setting, we will focus on examples of parabolic problems in order to employ a common framework for establishing well-definedness of the parameter-to-state map, the all-at-once version of the tangential condition trivially carries over to the wave equation (or also fractional sub- or superdiffusion) context by just replacing the first time derivative by a second (or fractional) one.

\medskip

The remainder of this paper is organized as follows. Section \ref{sec:aao} provides results for the all-at-once setting, that are also made use of in the subsequent Section \ref{Reduced} for the reduced setting. The proofs of the propositions in Section \ref{sec:aao} and the notation can be found in the appendix.


\section{All-at-once setting}\label{sec:aao}
The tangential cone condition and boundedness of the derivative in the all-at-once setting $\bF(\theta,u)=\by$ \eqref{bFthetauy} with 
\begin{equation}\label{bFaao}
\bF:\mathcal{X} \times \mathcal{U} \rightarrow \mathcal{W}^*\times H\times \mathcal{Y}
\,, \quad \bF(\theta,u)=\left(\begin{array}{c}\dot{u}-f(\theta,u)\\ u(0)-u_0\\ \calC(u)\end{array}\right)
\end{equation}
and the norms 
\[
\begin{aligned}
&\|(\theta,u)\|_{\calX\times\calU}:=\left(\|\theta\|_\calX^2+\|u\|_\calX^2\right)^{1/2}\,, \\
&\|(w,h,y)\|_{\calW^*\times H\times\calY}:=\left(\|w\|_{\calW^*}^2+\|h\|^2+\|y\|_\calY^2\right)^{1/2}\,, 
\end{aligned}
\]
on the product spaces
read as 
\begin{equation}\label{tangcone_aao_gen}
\begin{aligned}
&\|f(\theta,u)-f(\ttheta ,\tu )-f_\theta'(\theta,u)(\theta-\ttheta)-f_u'(\theta,u)(u-\tu)\|_{\mathcal{W}^*}\\
&\leq c_{tcc}^{AAO} \Bigl(
\|
{
\dot{u}-\dot{\tu}-f(\theta,u)+f(\ttheta ,\tu )
}
\|_{\mathcal{W}^*}^2 
+\|u(0)-\tu(0)\|_H^2
+\|\mathcal{C}(u-\tu)\|_{\mathcal{Y}}^2\Bigr)^{1/2}\,, \\
&\quad\forall (\theta,u),(\ttheta ,\tu ) \in B_\rho^{\mathcal{X}\times \mathcal{U}}(\theta^0,u^0)\,,
\end{aligned}
\end{equation}
and 
\begin{equation}\label{bFpbdd_aao_gen}
\begin{aligned}
&
\Bigl(\|\dot{v}-f_\theta'(\theta,u)\chi-f_u'(\theta,u)v\|_{\mathcal{W}^*}^2
+\|v(0)\|_H^2
+ \|\calC v\|_Y^2\Bigr)^{1/2}\\ 
&\leq C_{\bF} 
{
\Bigl(\|\chi\|_\calX^2+\|v\|_\calU^2\Bigr)^{1/2} 
}\,, \\
&\quad\forall (\theta,u)\in B_\rho^{\mathcal{X}\times \mathcal{U}}(\theta^0,u^0)\,,\ \chi\in \calX\,,\ v\in\calU
\end{aligned}
\end{equation}
where we have assumed linearity of $\mathcal{C}$.

Since the right hand side terms $\|u(0)-\tu(0)\|_H$ and $\|f(\theta,u)-f(\ttheta ,\tu )\|_{\mathcal{W}^*}$ in \eqref{tangcone_aao_gen} are usually too weak to help for verification of this condition, we will just skip it in the following
and consider 
\begin{equation}\label{tangcone_aao}
\begin{aligned}
&\|f(\theta,u)-f(\ttheta ,\tu )-f_\theta'(\theta,u)(\theta-\ttheta)-f_u'(\theta,u)(u-\tu)\|_{\mathcal{W}^*}\\
&\leq c_{tcc}^{AAO} \|\mathcal{C}(u-\tu)\|_{\mathcal{Y}}\,, \quad\forall (\theta,u),(\ttheta ,\tu ) \in B_\rho^{\mathcal{X}\times \mathcal{U}}(\theta^0,u^0)
\end{aligned}
\end{equation}
which under these conditions is obviously sufficient for \eqref{tangcone_aao_gen}.
Moreover, in order for the remaining right hand side term to be sufficiently strong 
in order to be able to dominate the left hand side,
we will need to have full observations in the sense that
\begin{equation}\label{fullobs}
\overline{\mathcal{R}(\mathcal{C}(t))}=Y\,.
\end{equation}
In the next section, it will be shown that under certain stability conditions on the generalized ODE in \eqref{model-1}, together with \eqref{fullobs}, the version \eqref{tangcone_aao} of the all-at-once tangential cone condition is sufficient for its reduced counterpart \eqref{tcc}.

Likewise, we will further consider the following sufficient conditions for boundedness of the derivative, 
\begin{equation}\label{bFpbdd_aao}
\begin{aligned}
&
\|f_\theta'(\theta,u)\|_{L(\calX,\mathcal{W}^*)}\leq C_{\bF,1}\,, \quad
&&\|f_u'(\theta,u)\|_{L(\calU,\mathcal{W}^*)}\leq C_{\bF,2}\,,\\
&\|\partial_t\|_{L(\calU,\mathcal{W}^*)}\leq C_{\bF,0}\,, \quad 
&&\|\calC\|_{L(\calU,
{
\calY
}
)}\leq C_{\bF,3}
\\
&\quad\forall (\theta,u)\in B_\rho^{\mathcal{X}\times \mathcal{U}}(\theta^0,u^0)\,.
\end{aligned}
\end{equation}

The function space setting considered here will be 
\begin{equation}\label{Uaao}
\begin{aligned}
&
\calU=\{u \in L^2(0,T;V): \dot{u} \in L^2(0,T;W^*)\}\hookrightarrow C(0,T;H)\,, \\
&\mathcal{W}=L^2(0,T;W)\,, \qquad
\mathcal{Y}=L^2(0,T;Y)\,, \ 
\end{aligned}
\end{equation}
so that the third bound in \eqref{bFpbdd_aao} is automatically satisfied with $C_{\bF,0}=1$.
We focus on Lebesgue and Sobolev spaces\footnote{In place of $V$, its intersection with $H_0^1(\Omega)$ might be considered in order to take into account homogeneous Dirichlet boundary conditions. For the estimates themselves, this does not change anything.}
\begin{equation}\label{VWY}
\begin{aligned}
&V=W^{s,m}(\Omega)\,, \ 
&&W=W^{t,n}(\Omega)\,, \ 
&&Y=L^q(\Omega)\,, 
\end{aligned}
\end{equation}
with $s,t\in[0,\infty)$, $m,n\in[1,\infty]$, $q\in[1,\hat{q}]$, and $\hat{q}$ the maximal index such that $V$ continuously embeds into $L^{\hat{q}}(\Omega)$, i.e. such that
\begin{equation}\label{qhat}
s-\frac{d}{m}\succeq -\frac{d}{\hat{q}}\,,
\end{equation}
so that with $\calC$ defined by the embedding operator $\calU\to\calY$, the last bound in \eqref{bFpbdd_aao} is automatically satisfied
\footnote{One could possibly think of also extending to more general Lebesgue spaces instead of $L^2$ with respect to time. As long as the summability index is the same for $\mathcal{W}$ and $\mathcal{Y}$ this would not change anything in Subsection \ref{subsec:bilin}. As soon as the summability indices differ, one has to think of continuity of the embedding $\mathcal{U}=L^{r_1}(0,T;V)\cap W^{1,r_2}(0,T;W^*)\hookrightarrow\mathcal{Y}=L^{r_3}(0,T;Y)$ as a whole, possibly taking advantage of some interpolation between $L^{r_1}(0,T;V)$ and $W^{1,r_2}(0,T;W^*)$. This could become very technical but might pay off in specific applications.}.
For the notation $\succeq$ we refer to the appendix.

The parameter space $\mathcal{X}$ may be very general at the beginning of Subsection \ref{subsec:bilin} and in Subsection \ref{subsec:invsource}. We will only specify it in the particular examples of Subsection \ref{subsec:bilin}.

We will now verify the conditions \eqref{tangcone_aao}, \eqref{bFpbdd_aao} for some (classes of) examples.

\subsection{Bilinear problems}\label{subsec:bilin}

Many coefficient identification problems in linear PDEs, such as the identification of a potential or of a diffusion coefficient, as mentioned above, can be treated in a general bilinear context.

Consider an evolution driven by a bilinear operator, i.e.,
\begin{equation}\label{bilin}
f(\theta,u)(t)=L(t)u(t)+((B\theta)(t))u(t)-g(t)\,,
\end{equation}
where for almost all $t\in(0,T)$, and all $\theta\in\mathcal{X}$, $v\in V$ we have $L(t),\,(B\theta)(t)\, \in \cL(V,W^*)$, $\theta\mapsto(B\theta)(t)v\,\in \cL(\calX,W^*)$, and $g(t)\in W^*$,
with 
\begin{equation}\label{CBCL}
\sup_{t\in[0,T]}\|L(t)\|_{\cL(V,W^*)}\leq C_L\,, \quad
\sup_{t\in[0,T]}\|(B\theta)(t)\|_{\cL(V,W^*)}\leq C_B\|\theta\|_\calX
\end{equation}
so that the first and second bounds in \eqref{bFpbdd_aao} are satisfied, due to the estimates
\[
\begin{aligned}
\|f_\theta'(\theta,u)\chi\|_{\calW^*}&=\left(\int_0^T\|((B
{
\chi
}
)(t))u(t)\|_{W^*}^2\right)^{1/2}\leq 
C_B \|\chi\|_\calX \left(\int_0^T\|u(t)\|_{V}^2\right)^{1/2}\\
\|f_u'(\theta,u)v\|_{\calW^*}&=\left(\int_0^T\|L(t)v(t)+((B\theta)(t))v(t)\|_{W^*}^2\right)^{1/2}\\
&\leq (C_L+C_B \|\theta\|_\calX) \left(\int_0^T\|v(t)\|_{V}^2\right)^{1/2}\,.
\end{aligned}
\]

For the left hand side in \eqref{tangcone_aao},
we have 
\[
\Bigl(f(\theta,u)-f(\ttheta,\tu)-f_u'(\theta,u)(u-\tu)-f_\theta'(\theta,u)(\theta-\ttheta)\Bigr)(t)=-((B(\theta-\ttheta)(t))(u-\tu)(t)\,,
\]
and \eqref{tangcone_aao} is satisfied if and only if 
\mpar{{removed ref. to \eqref{fullobs}}}
\[
\|(B(\theta-\ttheta))(u-\tu)\|_{\mathcal{W}^*}\leq c_{tcc}^{AAO}\|\mathcal{C}(u-\tu)\|_{\mathcal{Y}}
\,, \quad\forall (\theta,u),(\ttheta ,\tu ) \in B_\rho^{\mathcal{X}\times \mathcal{U}}(\theta^0,u^0)\,, 
\]
hold.
A sufficient condition for this to hold is
\begin{equation}\label{tangcone_aao_bilin}
\begin{aligned}
&\|(B(\theta-\ttheta))(t)(v-\tv)\|_{W^*}\leq c_{tcc}^{AAO}\|\mathcal{C}(t)(v-\tv)\|_Y\,,\\
& \qquad\forall (\theta,v),(\ttheta ,\tv ) \in B_\rho^{\mathcal{X}\times V}(\theta^0,u^0(t))\,,
\quad t\in(0,T)
\end{aligned}
\end{equation}

The proofs of the propositions for the following examples can be found in the appendix.
Likewise, the conditions on the summability and smoothness indices $s,t,p,q,m,n$ of the used spaces, 
\eqref{stmnpq_cprob}, \eqref{stmnpq_aprob}, 
\eqref{casegammkappage2}, \eqref{casegammkappale2}, \eqref{condts}, \eqref{condgrad}, 
\eqref{U1}, \eqref{U2}, \eqref{condPhi}, \eqref{condPsi}
as appearing in the formulation of the propositions, are derived there. 

\subsubsection*{Identification of a potential $c$}
Problem (\ref{c-problem-1})-(\ref{c-problem-3}) can be cast into the form \eqref{bilin} by setting $\theta=c$ and
\begin{equation}\label{LB_c-problem}
L(t)=\Delta\,, \quad (Bc)(t)v=
{
-cv
}, 
\end{equation}
(i.e., $(Bc)(t)$ is a multiplication operator with the multiplier $c$). We set
\begin{equation}\label{Xcprob}
\mathcal{X}=L^p(\Omega)\,.
\end{equation}

\begin{proposition}\label{prop:cprob_aao}
For $\calU$, $\calW$, $\calY$ according to \eqref{Uaao} with \eqref{VWY} \eqref{stmnpq_cprob}, $-\Delta\in\mathcal{L}(V,W^*)$, the operator $\bF$ defined by \eqref{bFaao}, \eqref{bilin}, \eqref{LB_c-problem}, $\calC=\mbox{id}:\calU\to \calY$ satisfies the tangential cone condition \eqref{tangcone_aao} with a uniformly bounded operator $\bF'(c)$,
i.e., the family of linear operators $(\bF'(c))_{c\in 
{
M
}
}$ is uniformly bounded in the operator norm, for $c$ in a bounded subset 
{
$M$
}
of $\calX$.

\end{proposition}

\begin{remark}
A full Hilbert space setting can be achieved by setting $p=q=m=n=2$ and choosing $s\geq0$, $t>\frac{d}{2}$.
\end{remark}

\subsubsection*{Identification of a diffusion coefficient $a$}
The $a$ problem  (\ref{a-problem-1})-(\ref{a-problem-3}) is defined by setting 
\begin{equation}\label{LB_a-problem}
L(t)\equiv0\,, \quad (Ba)(t)v=\nabla\cdot(a\nabla v), 
\end{equation}
so that 
\[
\begin{aligned}
&\|(B(\hat{a})(t))\hat{v}\|_{W^*}\\
&=\sup_{w\in W\,, \ \|w\|_W\leq1}\int_\Omega \hat{a}\nabla\hat{v}\cdot\nabla w\, dx
=\sup_{w\in W\,, \ \|w\|_W\leq1}\int_\Omega \hat{v}\Bigl(\nabla\hat{a}\cdot\nabla w+\hat{a}\Delta w\Bigr)\, dx\\
&\leq 
\|\hat{v}\|_{L^q} \Bigl(
\|\nabla\hat{a}\|_{L^p} \sup_{w\in W\,, \ \|w\|_W\leq1} \|\nabla w\|_{L^{\frac{p^* q}{q-p^*}}}
+\|\hat{a}\|_{L^r} \sup_{w\in W\,, \ \|w\|_W\leq1} \|\Delta w\|_{L^{\frac{r^* q}{q-r^*}}}\Bigr)\,.
\end{aligned}
\]
Note that since $Y=L^q(\Omega)$ we had to move all derviatives away from $\hat{v}$ by means of integration by parts, which forces us to use spaces of differentiability order at least two in $W$ and at least one in $\mathcal{X}$. Thus we here consider 
\begin{equation}\label{Xaprob}
\mathcal{X}=W^{1,p}(\Omega)\,.
\end{equation}

\begin{proposition}\label{prop:aprob_aao}
For $\calU$, $\calW$, $\calY$ according to \eqref{Uaao} with \eqref{VWY}, \eqref{stmnpq_aprob}, the operator $\bF$ defined by \eqref{bFaao}, \eqref{bilin}, \eqref{LB_a-problem}, $\calC=\mbox{id}:\calU\to \calY$ satisfies the tangential cone condition \eqref{tangcone_aao} with a uniformly bounded operator $\bF'(a)$.
\end{proposition}

\begin{remark}
A full Hilbert space setting $p=q=m=n=2$ requires to choose $s\geq0$ and $t\begin{cases}\geq 2\mbox{ if }d=1\\> 2\mbox{ if }d=2\\>1+\frac{d}{2}\mbox{ if }d\geq3\end{cases}$.
\end{remark}

\subsection{Nonlinear inverse source problems}\label{subsec:invsource}
Consider nonlinear evolutions that are linear with respect to the parameter $\theta$, i.e.
\begin{equation}\label{nonlin}
f(\theta,u)(t)=L(t)u(t)+\Phi(u(t))+\Psi(\nabla u(t))-B(t)\theta
\end{equation}
where for almost all $t\in(0,T)$, $L(t)\in \cL(V,W^*)$, $B(t)\in \cL(\mathcal{X},W^*)$ and $\Phi,\Psi\in C^2(\R)$ satisfy the H\"older continuity and growth conditions
\begin{equation}\label{growthPhipp}
|\Phi'(\lambda)-\Phi'(\tilde{\lambda})|\leq C_{\Phi''}(1+|\lambda|^{\gamma}+|\tilde{\lambda}|^{\gamma})|\tilde{\lambda}-\lambda|^\kappa
\end{equation}
for all $\tilde{\lambda},\lambda\in\R$
\begin{equation}\label{growthPsipp}
|\Psi'(\lambda)-\Psi'(\tilde{\lambda})|\leq C_{\Psi''}(1+|\lambda|^{\hat{\gamma}}+|\tilde{\lambda}|^{\hat{\gamma}})|\tilde{\lambda}-\lambda|^{\hat{\kappa}}
\end{equation}
for all $\tilde{\lambda},\lambda\in\R^d$, where $\gamma,\hat{\gamma},\kappa,\hat{\kappa}\geq0$.
We will show that the exponents $\gamma,\hat{\gamma}$ may actually be arbitrary as long as the smoothness $s,t$ of $V$ and $W$ is chosen appropriately.

\medskip

\begin{proposition}\label{prop:nlinvsource_aao}
The operator $\bF$ defined by \eqref{bFaao}, \eqref{bilin}, \eqref{LB_a-problem}, $\calC=\mbox{id}:\calU\to \calY$ in either of the four following cases
\begin{enumerate}
\item[(a)]
\eqref{growthPhipp} and $\Psi$ affinely linear and $\calU$, $\calW$, $\calY$ as in \eqref{Uaao} with \eqref{VWY}, 
\eqref{casegammkappage2}, \eqref{casegammkappale2}; 
\item[(b)]
\eqref{growthPhipp}, \eqref{growthPsipp} and $\calU$, $\calW$, $\calY$ as in \eqref{Uaao} with \eqref{VWY}, 
\eqref{casegammkappage2}, \eqref{casegammkappale2}, \eqref{condts}, \eqref{condgrad}; 
\item[(c)]
\eqref{growthPhipp} and $\Psi$ affinely linear, $\calW$, $\calY$ as in \eqref{Uaao}, $\calU$ as in \eqref{U1} with \eqref{VWY}, \eqref{qhat}, \eqref{condPhi}; 
\item[(d)]
\eqref{growthPhipp}, \eqref{growthPsipp}, $\calW$, $\calY$ as in \eqref{Uaao}, $\calU$ as in \eqref{U2} with \eqref{VWY}, \eqref{qhat}, \eqref{condPhi}, \eqref{condPsi}; 
\end{enumerate}
satisfies the tangential cone condition \eqref{tangcone_aao} with a uniformly bounded operator $\bF'(\theta)$.
\end{proposition}

\begin{remark}
A Hilbert space setting $p=q=m=n=2$ is therefore possible for arbitrary $\gamma,\kappa,\hat{\gamma}$, provided $t$ and $s$ are chosen sufficiently large, cf. \eqref{casegammkappage2}, \eqref{casegammkappale2} in case $C_{\Psi''}=0$, and additionally \eqref{condts}, \eqref{condgrad} otherwise.
\end{remark}

\section{Reduced setting} \label{Reduced}

In this section, we formulate the system (\ref{model-1})-(\ref{model-3}) by one operator mapping from the parameter space to the observation space. To this end, we introduce the parameter-to-state map
\[S:\calD\subseteq\calX\rightarrow\tilde{\calU}, \qquad\text{where}\qquad u=S(\theta) \text{ solves } (\ref{model-1})-(\ref{model-2})\]
then, with $\calD(F)=\calD$ the forward operator for the reduced setting can be expressed as
\begin{align} 
F:\calD(F)\subseteq\calX \rightarrow \calY, \qquad \theta \mapsto \calC(S(\theta)) 
\end{align}
and the inverse problem of recovering $\theta$ from $y$ can be written as
\begin{align*}
F(\theta)=y. \qquad
\end{align*}
Here, differently from the state space $\calU$ in the all-at-once setting, cf., \eqref{Uaao}, we use a non Hilbert state space  
\begin{align} \label{statespaceRe}
\tilde{\calU}=\{u \in L^\infty(0,T;V): \dot{u} \in L^2(0,T;W^*)\}
\end{align}
as this appears to be more appropriate for applying parabolic theory.
\bigbreak

We now establish a framework for verifying the tangential cone condition as well as boundedness of the derivative in this general setting.

For this purpose, we make the following assumptions.
\begin{assumption}   \label{assumption-R} \text{ }
\begin{enumerate}[label=(R\arabic*)]
\item \label{R1}
\textit{Local Lipschitz continuity of $f$}
\begin{align*}
\forall M\geq 0,& \exists L(M) \geq 0, \forall^{a.e.} t \in (0,T): \nonumber\\
& \|f(t,\theta_1,v_1)-f(t,\theta_2,v_2)\|_{W^*}
\leq L(M) (\|v_1-v_2\|_V +\|\theta_1-\theta_2\|_\calX),	\\
&\qquad\qquad\forall v_i \in V, \theta_i \in \calX: \|v_i\|_V, \|\theta_i\|_\calX \leq M, i=1,2. \nonumber
\end{align*}
\item \label{R2}
\textit{Well-definedness of the parameter-to-state map}
\[S:\calD(F)\subseteq\calX\rightarrow\tilde{\calU} \]
with $\tilde{\calU}$ as in (\ref{statespaceRe}) as well as its boundedness in the sense that there exists $C_S>0$ such that for all $\theta\in \calB_\rho^\calX(\theta^0)$ the estimate
\[
\|S(\theta)\|_{L^\infty(0,T;V)}\leq C_S
\]
holds.
\item \label{R3}
\textit{Continuous dependence on data of the solution to the linearized problem with zero initial data}, i.e., 
there exists a constant $C_{lin}$ such that for all $\theta\in\calB_\rho^\calX(\theta^0)$, $b\in \calW^*$, and any $z$ solving 
\begin{flalign}
& \dot{z}(t) = f'_u(\theta,S(\theta))(t)z(t)+b(t) \quad t\in (0,T) \qquad\qquad \label{linear-eq-1}\\ 
& z(0)=0,  \label{linear-eq-2}
\end{flalign}
the estimate
\begin{align} \label{linear-eq-estimate}
\|z\|_\calY \leq C_{lin}\|b\|_{\calW^*}. \qquad\qquad
\end{align}
holds.
\item \label{R4}
\textit{Tangential cone condition of the all-at-once setting \eqref{tangcone_aao}}
\begin{align*}
&\exists \rho >0, \forall 
{
(\theta,u),
}
(\tilde{\theta},\tilde{u})\in\calB_\rho^{\calX,\calU}(\theta^0,u^0): \\
&\|f(\tilde{\theta},\tilde{u})-f(\theta,u)-f'_u(\theta,u)(\tilde{u}-u)-f'_\theta(\theta,u)(\tilde{\theta}-\theta)\|_{\calW^*}
\leq c_{tcc}^{AAO} \|\calC\tilde{u}-\calC u\|_{\mathcal{Y}}. \qquad\qquad\quad
\end{align*}
\end{enumerate}
\end{assumption}\bigbreak

The main result of this section is as follows.

\begin{theo}
Suppose Assumption \ref{assumption-R} holds and $\calC$ is the embedding $V\hookrightarrow Y$.\\
Then there exists a constant $\rho>0$ such that for all $\theta, \tilde{\theta}\in\calB^\calX_\rho(\theta_0)\subset\calD(F)$,
\begin{enumerate}[label=\roman*)]
\item $F'(\theta)$ is uniformly bounded:
\[\|F'(\theta)\|_{
{
\mathcal{L}(\calX,\calY)
}
}
\leq M \]
\end{enumerate}
\quad for some constant $M$, and
\begin{enumerate}[label=ii)]
\item The tangential cone condition is satisfied:
\begin{align} \label{TCC-Re}
\|F(\tilde{\theta})-F(\theta)-F'(\theta)(\tilde{\theta}-\theta)\|_{\mathcal{Y}}
\leq c_{tcc}^{Re} \|F(\tilde{\theta})-F(\theta)\|_{\mathcal{Y}} 
\end{align}
for some small constant $c^{Re}_{tcc}$.
\end{enumerate}
\end{theo}

This is a consequence of the following two propositions, in which we combine the all-at-once versions of the tangential cone and boundedness conditons, respectively, with the assumed stability of $S$ and its linearization. 

\begin{proposition}\label{prop:tangcone}
Given $\calC$ is the embedding $V\hookrightarrow Y$ and $u_0$ is independent of $\theta$, the tangential cone condition in the reduced setting (\ref{TCC-Re}) follows from the one in the all-at-once setting \ref{R4} if the linearized forward operator is boundedly invertible as in \ref{R3} and $S$ is well defined according to \ref{R2}.
\end{proposition}
\proof

We begin by observing that the functions
\begin{flalign*}
 v&:= S(\theta)-S(\tilde{\theta})&\\
  w&:= S'(\theta)h&\\
 z&:= S(\theta)-S(\tilde{\theta})-S'(\theta)(\theta-\tilde{\theta})&
\end{flalign*}
solve the corresponding equations
\begin{alignat}{3}
\dot{v}(t)
&=f(\theta,S(\theta))(t) - f(\tilde{\theta},S(\tilde{\theta}))(t)  \qquad t\in(0,T), \qquad v(0)=0 \label{diffS}\\
 \dot{w}(t)
&=f'_u(\theta,S(\theta))w(t) + f'_\theta(\theta,S(\theta))h(t) \qquad t\in(0,T), \qquad w(0)=0 \label{S'}\\
\dot{z}(t)
&= f'_u(\theta,S(\theta))z(t) \nonumber\\
&\qquad +\big(-f'_u(\theta,S(\theta))v(t)- f'_\theta(\theta,S(\theta))(\theta-\tilde{\theta})(t) \\ 
&\qquad+   f(\theta,S(\theta))(t) -f(\tilde{\theta},S(\tilde{\theta}))(t) \big) \nonumber \qquad\qquad\\
&=: f'_u(\theta,S(\theta))z(t) + r(t) \qquad t\in(0,T), \qquad z(0)=0. \label{linearized-problem}
\end{alignat}
Hence we end up with the following estimate, using the assumed bounded invertibility of the linearized problem (\ref{linearized-problem}) and the fact that $\calC$ is the embedding $V\hookrightarrow Y$,
\begin{align}
\|F(\theta)-F(\tilde{\theta})-F'(\theta)(\theta-\tilde{\theta})\|_\calY
& =\|S(\theta)-S(\tilde{\theta})-S'(\theta)(\theta-\tilde{\theta})\|_\calY\nonumber\\
& \leq C_{lin} \|r\|_{\calW^*} \label{S-energy}\\
& \leq C_{lin} c_{tcc}^{AAO}\|F(\theta)-F(\tilde{\theta})\|_\calY, \label{TCC-AAO}
\end{align}
where $\|r\|_{\calW^*}$ and $c_{tcc}^{AAO}$ are respectively the left hand side and the constant in the all-at-once tangential cone estimate, applied to $u=S(\theta)$ and $\tilde{u}=S(\tilde{\theta})$.  \QED\\

\begin{remark}
The inverse problem \eqref{Fthetay} with \eqref{AC}, \eqref{S}, \eqref{FCS} can be written as a composition of the linear observation operator $\calC$ and the nonlinear parameter-to-state map $S$. 
Such problems have been considered and analyzed in \cite{Hof94}, but as opposed to that the inversion of our observation operator is ill-posed so the theory of \cite{Hof94} does not apply here.
\end{remark}

Note that in (\ref{TCC-AAO}), $c_{tcc}^{AAO}$ must be sufficiently small such that the tangential cone constant in the reduced setting $c_{tcc}^{Re}:=C_{lin} c_{tcc}^{AAO}$ fulfills the smallness condition required in convergence proofs as well.
Moreover we wish to emphasize that for the proof of Proposition \ref{prop:tangcone}, the constant $C_{lin}$ does not need to be uniform but could as well depend on $\theta$.
Also the uniform boundedness condition on $S$ from \ref{R2} is not yet needed here.

Under further assumptions on the defining functions $f$, we also get existence and uniform boundedness of the linear operator $F'(\theta)$ as follows.
\begin{proposition}
Let $S$ be well defined and bounded according to \ref{R2}, and let \ref{R1}, \ref{R3} be satisfied.\\
Then $F'(\theta)$ is G\^ateaux differentiable and its derivative given by
\begin{equation}\label{F'}
F'(\theta) : \calX \rightarrow \calY,
\qquad\text{where}\qquad F'(\theta)h=w  \text{ solves } (\ref{S'}) 
\end{equation}
is uniformly bounded in $\calB_\rho^\calX(\theta_0)$.
\end{proposition}
\proof
For differentiablity  of $F$ relying on conditions \ref{R1}-\ref{R3}, we refer to \cite[Proposition 4.2]{Nguyen}.
Moreover again using \ref{R1}-\ref{R3}, for any $\theta\in\calB_\rho^\calX(\theta_0)$ we get
\begin{align*}
\|F'(\theta)h\|_\calY 
&=\|S'(\theta)h\|_\calY \leq C_{lin}\|f'_\theta(\theta,S(\theta))h\|_{L^2(0,T;W^*)} \\
&\leq C_{lin}\sqrt{T}\|f'_\theta(\theta,S(\theta))\|_{\calX\rightarrow W^*}\|h\|_\calX\\
&\leq C_{lin}\sqrt{T}L(M) \|h\|_\calX
\end{align*}
for $M=C_S+\|\theta_0\|_\calX+\rho$, where $L(M)$ is the Lipschitz constant in \ref{R1} and $C_{lin}$ is as in \ref{R3}. Above, we employ boundedness of $S$ by $C_S$ as assumed in \ref{R2}.\\
This proves uniform boundedness of $F'(\theta)$.\\

We now discuss Assumption \ref{assumption-R} in more detail.

\begin{remark} \label{VeqW}
For the case $V=W$.\\
We rely on the setting of a Gelfand triple $V\subseteq H\subseteq V^*$ for the general framework of nonlinear evolution equations.
By this, \ref{R2} can be fulfilled under the conditions suggested by Roub\' i\v cek \cite[Theorems 8.27, 8.31]{Roubicek}: \medbreak
For every $\theta \in\calD(F)$
\begin{enumerate}[label=(S\arabic*)]
\item \label{S1}
and for almost $t \in (0,T)$, the mapping $-f(t,\theta,\cdot)$ is pseudomonotone, {i.e.,
$-f(t,\theta,\cdot)$ is bounded and
\begin{align*}
\left.\begin{aligned}
\underset{k\rightarrow\infty}{\liminf} \dupair{f(t,\theta,u_k),u_k-u} \geq 0 \\
\qquad\qquad\qquad\qquad\quad u_k \rightharpoonup u
\end{aligned}\right\}
\Rightarrow
\begin{cases}
\dupair{f(t,\theta,u),u-v} \geq \underset{k\rightarrow\infty}{\limsup}\dupair{f(t,\theta,u_k),u_k-v}\\
\forall v \in V.
\end{cases}
\end{align*}
}
\item \label{S2}
$ -f(\cdot,\theta,\cdot)$ is semi-coercive, i.e.,
\[ \forall v \in V, \forall^{a.e.}t \in (0,T): \dupair{ -f(t,\theta,v),v}_{V^*,V} \geq C_0^\theta|v|^2_V -C_1^\theta(t)|v|_V-C_2^\theta(t)\|v\|_H^2 \qquad\]
for some $C_0^\theta > 0, C_1^\theta \in L^2(0,T), C_2^\theta \in L^1(0,T)$ and some seminorm $|.|_V$ satisfying\\ $\forall v \in V: \|v\|_V \leq c_{|.|}(|v|_V + \|v\|_H)$ for some $c_{|.|}>0.$
\item \label{S3}
$f$ satisfies {the growth condition 
\begin{align*}
\exists\gamma^\theta\in L^2(0,T), \hbar^\theta:\R\rightarrow\R \text{ increasing}: \|	f(t,\theta,v)\|_{V^*}\leq\hbar^\theta(\|v\|_H)(\gamma^\theta(t)+\|v\|_V)
\end{align*}}
and a condition for uniqueness of the solution, e.g.,
\[ \forall u,v \in V, \forall^{a.e.}t \in (0,T): \langle f(t,\theta,u) - f(t,\theta,v) ,u - v\rangle_{V^*,V} \leq \rho^\theta(t)\|u - v\|^2_H \qquad\]
for some $\rho^\theta \in L^1(0,T)$
\end{enumerate}
and further conditions for $S(\theta)\in L^\infty(0,T;V)$, e.g, \cite[Theorem 8.16, 8.18]{Roubicek}.

In case of linear $f(t,\theta,\cdot)$, \ref{S1}-\ref{S3} boil down to boundedness and semi-coercivity \ref{S2} of $-f(\cdot,\theta,\cdot)$ according to \cite[Theorem 8.27, 8.31, 8.28]{Roubicek}. Alternatively, one can observe that linear boundedness implies the growth condition in \ref{S3}  with $ \gamma^\theta=0, \hbar^\theta=\|f(\cdot,\theta)\|$, and \ref{S2} implies the rest of \ref{S3} with $\rho^\theta=C^\theta_2$ if $C^\theta_1\leq0$ as $C^\theta_0>0$. The pseudomonotonicity assumption \ref{S1}, which guarantees week convergence of $f(\cdot,\theta,u_k)$ to $f(\cdot,\theta,u)$ when the approximation solution sequence $u_k$ converges weakly to $u$, can  be replaced by weak continuity of $f(\cdot,\theta,\cdot)$ which holds in this linear bounded case.

Treating the linearized problem (\ref{linear-eq-1})-(\ref{linear-eq-2}) as an independent problem, we can impose on $f'_u(\theta,S(\theta))$ the {boundedness} and semi-coercivity properties, then \ref{R3}  follows. 
\end{remark}

\begin{remark} \label{VneqW}
For general spaces $V, W.$\\
Some examples even in case $V\not=W$ allow to use the results quoted in Remark \ref{VeqW} with an appropriately chosen Gelfand triple, see, e.g., Subsection \ref{c-problem} below.

When dealing with linear and quasilinear parabolic problems, detailed discussions for unique exsistence of the solution are exposed in the books, e.g., of Evans \cite{Evans}, Ladyzhenskaya et al. \cite{{LadyzhenskayaSolonnikovUraltseva}}, Pao \cite{Pao}. If constructing the solution to the initial value problem through the semigroup approach, one can find several results, e.g, from Evans \cite{Evans}, Pazy \cite{Pazy} combined with the elliptic results from Ladyzhenskaya et al. \cite{{LadyzhenskayaUraltseva}}.

Addressing \ref{R3}, a possible strategy is using the following dual argument.\\
Suppose $W$ is reflexive and $z$ is a solution to the problem (\ref{linear-eq-1})-(\ref{linear-eq-2}), then by the Hahn-Banach Theorem
\begin{flalign*}
\qquad \|z\|_{L^2(0,T;V)}
&=\sup_{\|\phi\|_{L^2(0,T;V^*)}\leq 1}\int_0^T \langle z, \phi \rangle_{V,V^*} dt & \\
&=\sup_{\|\phi\|_{L^2(0,T;V^*)}\leq 1}\int_0^T \langle z, -\dot{p} - f'_u(\theta,S(\theta))^*p \rangle_{V,V^*} dt &\\
&=\sup_{\|\phi\|_{L^2(0,T;V^*)}\leq 1}\int_0^T \langle\dot{z} - f'_u(\theta,S(\theta))z, p \rangle_{W^*,W} dt &\\
&=\sup_{\|\phi\|_{L^2(0,T;V^*)}\leq 1}\int_0^T \langle b, p \rangle_{W^*,W} dt &\\
&\leq \sup_{\|\phi\|_{L^2(0,T;V^*)}\leq 1} \|b\|_{L^2(0,T;W^*)} \|p\|_{L^2(0,T;W)},&
\end{flalign*}
where
\begin{flalign*}
f'_u(\theta,S(\theta))(t): V \rightarrow W^*,\qquad\qquad f'_u(\theta,S(\theta))(t)^*: W^{**}=W \rightarrow V^*,&
\end{flalign*}
and $p$ solves the adjoint equation
\begin{align}
-\dot{p}(t)& = f'_u(\theta,S(\theta))^*p(t)+\phi(t) \quad t\in (0,T) & \label{adjoin-eq-1}\\ 
 p(T)&=0. \label{adjoin-eq-2}&
\end{align}

If in the adjoint problem the estimate 
\begin{align} \label{estadjoint}
\|p\|_{L^2(0,T;W)}\leq \tilde{C}_{lin}\|\phi\|_{L^2(0,T;V^*)}
\end{align}
holds for some uniform constant $\tilde{C}_{lin}$, then we obtain
\begin{align} 
\|z\|_\calY\leq \|\calC\|_{V\to Y}\, \|z\|_{L^2(0,T;V)} \leq \|\calC\|_{V\to Y}\, \tilde{C}_{lin}\|b\|_{\calW^*}. 
\end{align}
Thus \ref{R3} is fulfilled.

\medskip

So we can replace \ref{R3} by
\begin{enumerate}[label=(R3-dual),leftmargin=1.6cm] 
\item \label{R3-dual}\textit{Continuous dependence on data of the solution to the adjoint linearized problem} associated with zero final condition, i.e., 
there exists a constant $\tilde{C}_{lin}$ such that for all $\theta\in\calB_\rho^\calX(\theta^0)$, $\phi\in L^2(0,T,V^*)$, and any $p$ solving \eqref{adjoin-eq-1}-\eqref{adjoin-eq-2}, the estimate \eqref{estadjoint} holds.
\end{enumerate}

\end{remark}

In the following sections, we examine the specific examples introduced in the introduction, in the relevant function space setting
\begin{align}
&\calX = L^p(\Omega) \quad\text{or}\quad \calX=W^{1,p}(\Omega)  \qquad\qquad\qquad p \in [1,\infty]\\ 
& Y = L^q(\Omega) \qquad\qquad\qquad\qquad\qquad\qquad\qquad\hphantom{l} q \in [1,\bar{q}] \qquad\\
&\tilde{\calU} = \{u \in L^\infty(0,T;V): \dot{u} \in L^2(0,T;W^*)\},
\end{align}
where $V, W$ will be chosen subject to the particular example, where $\hat{q}$ is the maximum power allowing $V\hookrightarrow L^{\hat{q}}(\Omega)$ and $\bar{q}\leq\hat{q}$ is the maximum power such that (\ref{linear-eq-estimate}) in \ref{R3} holds.

\subsection{Identification of a potential} \label{c-problem}

We investigate this  problem in the function spaces
\begin{align*} 
 \mathcal{D}(F)=\calX=L^p(\Omega),  \qquad Y=L^q(\Omega),  \qquad V=L^2(\Omega), \qquad W=H^2(\Omega)\cap H_0^1(\Omega). 
\end{align*}\bigbreak

Now we verify the conditions proposed in Assumption \ref{assumption-R}.

\begin{enumerate}[label=(R\arabic*)]
\item \textit{Local Lipschitz continuity of $f$:}\medbreak
Applying H\"older's inequality, we have
\begin{align*}
\|f(\tilde{c},\tilde{u})-f(c,u)\|_{W^*} &=\|\tilde{c}\tilde{u}-cu\|_{W^*} 
= \sup_{\|w\|_W \leq 1} \int_\Omega (\tilde{c}\tilde{u}-cu)w dx \\
&\leq \sup_{\|w\|_W \leq 1} \|w\|_W C_{W\rightarrow L^{\bar{p}}}\left( \int_\Omega | \tilde{c}(\tilde{u}-u)+ (\tilde{c}-c)u|^{\bar{p}^*}dx\right)^{\frac{1}{\bar{p}^*}} \qquad\qquad\qquad\\
& \leq C_{W\rightarrow L^{\bar{p}}}( \|\tilde{c}\|_{L^p}\|\tilde{u}-u\|_{L^r}+\|\tilde{c}-c\|_{L^p}\|u\|_{L^r} )\\
& \leq L(M)(\|\tilde{u}-u\|_V + \|\tilde{c}-c\|_\calX)
\end{align*}
with the dual index $\bar{p}^*=\frac{\bar{p}}{\bar{p}-1}$ and $r=\frac{\bar{p}p}{\bar{p}p-p-\bar{p}}, L(M)=C_{W\rightarrow L^{\bar{p}}} C_{V\rightarrow L^r}(\|u\|_V+\|\tilde{u}\|_V + \|c\|_\calX+\|\tilde{c}\|_\calX) + 1$. Above, we invoke the continuous embbedings through the constants $C_{W\rightarrow L^{\bar{p}}}, C_{V\rightarrow L^r}$, where $\bar{p}$ denotes the maximum power allowing $W\subseteq L^{\bar{p}}$. Thus we are supposing
\begin{align} \label{c-cond-1}
p\geq\max\left\{\frac{2\bar{p}}{\bar{p}-2}\,, \ \frac{\bar{p}}{\bar{p}-1}\right\}=\frac{2\bar{p}}{\bar{p}-2}\mbox{ and }
2-\frac{d}{2}\succeq-\frac{d}{\bar{p}} 
\end{align}
in order to guarantee $V=L^2(\Omega) \hookrightarrow L^r(\Omega)$ and $W=H^2(\Omega)\cap H_0^1(\Omega) \hookrightarrow L^{\bar{p}}(\Omega)$

\item \textit{Well-definedness and boundedness of the parameter-to-state map:}\medbreak

Verifying {boundedness and semi-coercivity} conditions
with the Gelfand triple $H_0^1(\Omega)\hookrightarrow L^2(\Omega)\hookrightarrow H^{-1}(\Omega)$ (while remaining with $V=L^2(\Omega)$ in the definition of the space $\tilde{\calU}$) 
shows that, for $u_0\in L^2(\Omega), \varphi\in L^2(0,T;H^{-1}(\Omega))$ 
the initial value problem (\ref{c-problem-1})-(\ref{c-problem-3}) admits a unique solution $u\in W(0,T):= \{u \in L^2(0,T;H_0^1(\Omega)): \dot{u} \in L^2(0,T;H^{-1}(\Omega)\}\subset \{u \in L^\infty(0,T;L^2(\Omega)): \dot{u} \in L^2(0,T;H^{-2}(\Omega))\}=\tilde{\calU}.$

\quad Indeed, coercivity is deduced as follows. For
\begin{flalign} \label{c-cond-2}
&\qquad\qquad\quad p\geq 2,\quad d\leq 3,&
\end{flalign}
we see
\begin{align} \label{c-coercive}
\int_\Omega cu^2 dx
&\leq \|c\|_{L^2(\Omega)} \left(\int_\Omega u^4dx\right)^{\frac{1}{2}}
\leq \|c\|_{L^2(\Omega)} \left(\int_\Omega u^2dx\right)^{\frac{1}{4}} \left(\int_\Omega u^6dx\right)^{\frac{1}{4}} \qquad\qquad\qquad\nonumber\\
&\leq \|c\|_{L^2(\Omega)} \|u\|_{L^2(\Omega)}^{\frac{1}{2}} \|u\|_{L^6(\Omega)}^{\frac{3}{2}}\nonumber\\
&\leq C_{H_0^1\rightarrow L^6}\|c\|_{L^2(\Omega)} \|u\|_{L^2(\Omega)}^{\frac{1}{2}} \|u\|_{H_0^1(\Omega)}^{\frac{3}{2}}\\
&\leq  C_{H_0^1\rightarrow L^6}\|c\|_{L^2(\Omega)}\left( \frac{1}{4\epsilon}\|u\|_{L^2(\Omega)}\|u\|_{H_0^1(\Omega)} + \epsilon\|u\|_{H_0^1(\Omega)}^2 \right)\nonumber\\
&\leq  C_{H_0^1\rightarrow L^6}\|c\|_{L^2(\Omega)}\left( \frac{1}{16\epsilon\epsilon_1}\|u\|^2_{L^2(\Omega)} + \frac{\epsilon_1}{4\epsilon}\|u\|_{H_0^1(\Omega)}^2 + \epsilon\|u\|_{H_0^1(\Omega)}^2 \right), \nonumber
\end{align}
which yields semi-coercivity
\begin{align*}
&\langle -f(t,c,u),u \rangle_{H^{-1},H_0^1}
=\int_\Omega (-\Delta u+cu)u dx\\
&\ \geq \left( 1- C_{H_0^1\rightarrow L^6}\|c\|_{L^2(\Omega)}\left(\frac{\epsilon_1}{4\epsilon}+\epsilon\right)\right)\|u\|_{H_0^1(\Omega)}^2 -  \frac{C_{H_0^1\rightarrow L^6}}{16\epsilon\epsilon_1}\|c\|_{L^2(\Omega)}\|u\|^2_{L^2(\Omega)}, \\
&=:C_0^c \|u\|_{H_0^1(\Omega)}^2 + C_2^c \|u\|^2_{L^2(\Omega)},
\end{align*}
where the constant $C_0^c$ is positive if choosing $\epsilon_1 <\epsilon$ and  $\epsilon, \epsilon_1$ sufficiently small.

Boundedness of $f$ can be concluded from
\begin{align*}
\|-f(&t,c,u) \|_{H^{-1}(\Omega)}=\sup_{\|v\|_{H_0^1}\leq 1}\int_\Omega (-\Delta u+cu)v dx\\
&\leq\sup_{\|v\|_{H_0^1}\leq 1}\left( \|u\|_{H^1(\Omega)}\|v\|_{H^1(\Omega)}+C_{H_0^1\rightarrow L^6}C_{H_0^1\rightarrow L^3}\|c\|_{L^2(\Omega)} \|u\|_{H^1(\Omega)}\|v\|_{H^1(\Omega)}\right)\\
&\leq C \|c\|_{L^2(\Omega)} \|u\|_{H^1(\Omega)}.
\end{align*}

Moreover, by the triangle inequality: $\|c\|_{L^2(\Omega)}\leq \|c^0\|_{L^2(\Omega)}+\|c-c^0\|_{L^2(\Omega)}\leq \|c^0\|_{L^2(\Omega)}+\rho$, semi-coercivity of $f$ is satisfied with the constants $C_0, C_1$ now depending only on the point $c^0$. This hence gives us uniform boundedness of $S$ on the ball $\calB^\calX_\rho(c^0)$.

\item \textit{Continuous dependence on data of the solution to the linearized problem with zero initial data:}\medbreak
We use the duality argument mentioned in Remark \ref{VneqW}. To do so, we need to prove existence of the adjoint state $p\in L^2(0,T;W)$ and the associated estimate \ref{R3-dual}.

\quad Initially, by the transformation $v=e^{-\lambda t}p$ and putting $\tau=T-t$, the adjoint problem (\ref{adjoin-eq-1})-(\ref{adjoin-eq-2}) is equivalent to
\begin{flalign} 
&\qquad\qquad\qquad\dot{v}(t)-\Delta v(t) +(\lambda+c)v(t)= e^{-\lambda t}\phi(t) \qquad t \in (0,T)& \label{c-problem-1'}\\
&\qquad\qquad\qquad v(0)= 0 . \label{c-problem-2'}
\end{flalign}
We note that this problem with $c=\hat{c} \in L^\infty(\Omega), \lambda+\hat{c}>-C_{PF}$, the constant in the Poincar\' e-Friedrichs inequality, $\phi\in L^2(0,T;L^2(\Omega)), \partial\Omega \in C^2$,
admits a unique solution in $L^2(0,T;H^2(\Omega)\cap H_0^1(\Omega))$ \cite[Section 7.1.3, Theorem 5]{Evans}
\footnote{ 
{
where smoothness of the domain can be slightly relaxed to $C^{1,1}$ as assumed here, see, e.g., \cite{Grisvard}
}
}
and the operator $\frac{d}{dt}-\Delta + (\lambda+\hat{c}): L^2(0,T;H^2(\Omega)\cap H_0^1(\Omega))\rightarrow L^2(0,T;L^2(\Omega))\times H^1(\Omega), p \mapsto (\phi,p_0)$ is boundedly invertible.

\quad Suppose $u$ solves (\ref{c-problem-1'})-(\ref{c-problem-2'}), by the identity
\begin{flalign*}
\qquad \dot{u}&-\Delta u  + (\lambda + c)u = e^{-\lambda t}\phi \quad \Leftrightarrow \quad \dot{u}-\Delta u  + (\lambda +\hat{c})u = e^{-\lambda t}\phi + (\hat{c}-c)u  &\\
\quad u &= \left(\frac{d}{dt}-\Delta + (\lambda+\hat{c})\right)^{-1}\left[e^{-\lambda t}\phi+(\hat{c}-c)u\right]&\\
&=: \T u,&
\end{flalign*}
we observe that  $\T: L^2(0,T;H^2(\Omega)\cap H_0^1(\Omega)) \rightarrow L^2(0,T;H^2(\Omega)\cap H_0^1(\Omega))$
is a contraction
\begin{align} \label{c-contraction}
&\|\T( u -v)\|_{L^2(0,T; H^2\cap H_0^1)} \nonumber\\
&\leq \left\|\left(\frac{d}{dt}-\Delta + (\lambda+\hat{c})\right)^{-1}\right\|_{L^2(0,T;L^2(\Omega))\rightarrow L^2(0,T;H^2\cap H_0^1)}
\hspace*{-2cm}\|(\hat{c}-c)(u-v)\|_{L^2(0,T;L^2(\Omega))}  \nonumber\\
&\leq C^{\hat{c}}\|\hat{c}-c\|_{L^p}\|u-v\|_{L^2(0,T;L^{\frac{2p}{p-2}}(\Omega))} \nonumber\\
&\leq C\epsilon \|u-v\|_{L^2(0,T;H^2\cap H_0^1)} ,
\end{align}
where $C\epsilon<1$ if we assume $\hat{c}=c^0\in L^\infty(\Omega)$ and $\rho$ is sufficiently small. In some case, smallness of $\rho$ can be omitted (discussed at the end of \ref{R3}). Estimate (\ref{c-contraction}) holds provided
\begin{flalign} \label{c-cond-3}
&\qquad\qquad\quad W= H^2(\Omega) \cap H_0^1(\Omega) \hookrightarrow L^{\frac{2p}{p-2}}(\Omega)\mbox{ i.e., } \quad p\geq\frac{2\bar{p}}{\bar{p}-2}.&
\end{flalign}
Thus, for $\phi\in L^2(0,T;L^2(\Omega))$ there exists a unique solution $v \in L^2(0,T;H^2(\Omega)\cap H_0^1(\Omega))$ to the problem (\ref{c-problem-1'})-(\ref{c-problem-2'}), which implies $p=e^{\lambda t}v \in L^2(0,T;H^2(\Omega) \cap H_0^1(\Omega))$ is the solution to the adjoint problem (\ref{adjoin-eq-1})-(\ref{adjoin-eq-2}).\\

\quad Observing that $p$ solves
\begin{flalign*} 
&\qquad\qquad\qquad\dot{p}(t)-\Delta p(t) + \hat{c}p(t) = (\hat{c}-c)p(t) + \phi(t) \qquad t \in (0,T)& \\
&\qquad\qquad\qquad p(0)=0, &	
\end{flalign*}
employing again \cite[Section 7.1.3 , Theorem 5]{Evans} and {smallness of $\rho$} yields
\begin{align} \label{p-energy}
 \|p\|_{L^2(0,T;W)}
&\leq C (\|(\hat{c}-c)p\|_{L^2(0,T;L^2(\Omega))} + \|\phi\|_{L^2(0,T;L^2(\Omega))}) \nonumber\\
&\leq C (2\rho\|p\|_{L^2(0,T;H^2\cap H_0^1)} + \|\phi\|_{L^2(0,T;L^2(\Omega))} )\\
&\leq C \|\phi\|_{L^2(0,T;V^*)} \nonumber
\end{align}
with some constant $C$ independent of $\theta\in \calB^\calX_\rho(c^0)$. This yields \ref{R3-dual} with $\bar{q}=2$.

If $d=1, p=2$ or $d=2, p>2$ or $d=3, p\geq\frac{12}{5}$, the smallness condition on $\rho$ can be omitted. Indeed, for $d=3, p\geq\frac{12}{5}$ testing the adjoint equation by $-\Delta p$ yields
\begin{align} \label{c-rho-1}
\int_\Omega-\dot{p}\Delta p &+ (\Delta p)^2dx
=\int_\Omega(cp-\phi)\Delta p dx \nonumber\\
\frac{1}{2}\frac{d}{dt}\|\nabla p\|^2_{L^2(\Omega)}&+\|\Delta p\|^2_{L^2(\Omega)}
\leq \frac{1}{2}\|\Delta p\|^2_{L^2(\Omega)} +\|\phi\|^2_{L^2(\Omega)}+\|cp\|^2_{L^2(\Omega)}\\
\frac{1}{2}\frac{d}{dt}\|\nabla p\|^2_{L^2(\Omega)}&+\frac{1}{2}\|\Delta p\|^2_{L^2(\Omega)}
\leq \|\phi\|^2_{L^2(\Omega)} + \|c\|^2_{L^p(\Omega)}\left(\int_\Omega p^{\frac{p}{p-2}+\frac{p}{p-2}}dx\right)^{\frac{p-2}{p}} \nonumber\\
&\leq \|\phi\|^2_{L^2(\Omega)} + \|c\|^2_{L^p(\Omega)}\|p\|_{L^6(\Omega)}\|p\|_{L^\infty(\Omega)}|\Omega|^{\frac{5p-12}{6p}} \nonumber\\
&\leq \|\phi\|^2_{L^2(\Omega)} + (\|c^0\|^2_\calX+\rho^2)|\Omega|^{\frac{5p-12}{6p}}\bigg(\frac{C^2_{H_0^1\rightarrow L^6}}{4\epsilon}\|\nabla p\|^2_{L^2(\Omega)}\nonumber\\
&\qquad+ 
\epsilon C^2_{H^2\cap H_0^1\rightarrow L^\infty}\left(\|\Delta p\|^2_{L^2(\Omega)}+\|\nabla p\|^2_{L^2(\Omega)}\right)\bigg), \nonumber
\end{align}
{where in the last estimate we apply Young's inequality.} Choosing $\epsilon$ sufficiently small allows us to subtract the term involving $\|\Delta p\|^2_{L^2(\Omega)}$ on the right hand side from the one on the left hand side and get a positive coefficient in front. Here, the choice of $\epsilon$ depends only on the constants $c^0, \rho, \Omega, C_{H^2\cap H_0^1\rightarrow L^\infty}$.\\
It is also obvious that, if $d<3$, in the second line of the above calculation, we can directly estimate as follow 
\begin{equation}
\begin{split}\label{c-rho-2}
&d=1, p=2: \|cp\|^2_{L^2(\Omega)}\leq\|c\|^2_{L^2(\Omega)}\|p\|^2_{L^\infty(\Omega)}\leq C^2_{H_0^1\rightarrow L^\infty}\|c\|^2_{L^2(\Omega)}\|\nabla p\|^2_{L^2(\Omega)}\\
&d=2, p>2:  \|cp\|^2_{L^2(\Omega)}\leq\|c\|^2_{L^p(\Omega)}\|p\|^2_{L^{\frac{2p}{p-2}}(\Omega)} \!\!\!\!\!\!\leq C^2_{H_0^1\rightarrow L^{\frac{2p}{p-2}}}\!\!\!\|c\|^2_{L^p(\Omega)}\|\nabla p\|^2_{L^2(\Omega)} .
\end{split}
\end{equation} 
Employing firstly Gronwall-Bellman inequality with initial data $\nabla p(0)=0$, then taking the integral on $[0,T]$, we obtain
\begin{align} \label{c-regularity}
\|p\|_{L^\infty(0,T;H^1(\Omega))}+\|\Delta p\|_{L^2(0,T;L^2(\Omega))}\leq C \|\phi\|_{L^2(0,T;L^2(\Omega))} 
\end{align}
with the constant $C$ depending only on $c^0, \rho$. This estimate is valid for all $c\in\calB^\calX_\rho(c^0)$.
Since the adjoint problem has the same form as the original problem, applying \eqref{c-regularity} in \eqref{c-contraction} we can relax $\hat{c}$, by means of without fixing $\hat{c}=c^0$ but chossing it sufficiently close to $c$ since $\overline{L^\infty(\Omega)}=L^p(\Omega), |\Omega|<\infty$ to have $C^{\hat{c}}\epsilon\leq C\epsilon$ arbitrarily small with constant $C$ as in \eqref{c-regularity}. Therefore the constraint on smallness of $\rho$ can be omitted in these cases. 

\item \textit{All-at-once tangential cone condition:} \medbreak
According to \eqref{qhat}, \eqref{stmnpq_cprob} with $s=0$, $t=2$, $m=n=2$, this follows if
\[
\frac{p}{p-1}\leq q\leq\hat{q}\geq2\mbox{ and }2-\frac{d}{2}\succeq-\frac{d(p-1)}{p}+\frac{d}{q}.
\]
\end{enumerate}

\begin{corollary}
Assume $u_0\in L^2(\Omega), \varphi\in L^2(0,T;H^{-1}(\Omega))$, and\\
\begin{align} \label{c-space}
\calD(F) =\calX=L^p(\Omega),\quad Y=L^q(\Omega),\quad &V=L^2(\Omega),\quad W=H^2(\Omega)\cap H_0^1(\Omega)\nonumber\\  &p\geq 2,\quad q\in \left[\underline{q},2\right],\quad d\leq 3 
\end{align}
with $\underline{q}=\max\left\{\frac{p}{p-1}\,, \ \min\limits_{q\in[1,\infty]} \left\{2-\frac{d}{2}\succeq -\frac{d(pq-p-q)}{pq}\right\}\right\}$.

Then $F$ defined by $F(c)=u$ solving \eqref{c-problem-1}-\eqref{c-problem-3} satisfies the tangential cone condition (\ref{TCC-Re}) with a uniformly bounded operator $F'(c)$ defined by (\ref{F'}), see also \cite{HNS95} for the static case.
\end{corollary}

\begin{remark}
This allows a full Hilbert space setting of $X$ and $Y$ by choosing $p=q=2$ as long as $d\leq 3$.
\end{remark}

\subsection{Identification of a diffusion coefficient} \label{a-problem}

We pose this  problem in the function spaces
\begin{align} \label{a-cond-1}
&\calX=W^{1,p}(\Omega),  \quad Y=L^q(\Omega),  \quad  V=L^2(\Omega), \quad W=H^2(\Omega)\cap H_0^1(\Omega)\qquad p>d &
\end{align} 
so that $\calX\hookrightarrow L^\infty(\Omega)$ and define the domain of $F$ by
\begin{align} \label{a-problem-D(F)}
\calD(F)=\{a \in \calX: a\geq\underline{a}>0 \text{ a.e. on } \Omega \}. \qquad
\end{align}\bigbreak

Now we examine the conditions \ref{R1}-\ref{R3}.

\begin{enumerate}[label=(R\arabic*)]
\item \textit{Local Lipschitz continuity of $f$:}
\begin{align*}
\|-\nabla\cdot &\Bigl(\tilde{a}\nabla \tilde{u}\Bigr) +\nabla\cdot\Bigl(a\nabla u\Bigr)\|_{W^*}\\ 
&= \sup_{\|w\|_W \leq 1} \int_\Omega (\tilde{a}\nabla \tilde{u} - a\nabla u)\nabla w dx\\
&= \sup_{\|w\|_W \leq 1} \int_\Omega ( a\nabla(\tilde{u}-u)+ (\tilde{a}-a)\nabla \tilde{u})\nabla w dx\qquad\qquad\qquad\qquad\qquad\\
&= \sup_{\|w\|_W \leq 1} \int_\Omega (\tilde{u}-u) (\nabla a\nabla w + a\Delta w)+ \tilde{u}(\nabla(\tilde{a}-a)\nabla w + (\tilde{a}-a)\Delta w) dx \qquad\\
&\leq \sup_{\|w\|_W \leq 1} \int_\Omega (\|\tilde{u}-u\|_{L^2} \|\nabla a\|_{L^p}+\|\tilde{u}\|_{L^2}\|\nabla(a-a)\|_{L^p})\|\nabla w\|_{L^{\frac{2p}{p-2}}}\\
&\qquad\qquad\qquad\quad + (\|\tilde{u}-u\|_{L^2}\|a\|_{L^\infty}+\|\tilde{u}\|_{L^2}\|a-a\|_{L^\infty} )\|\Delta w\|_{L^2}dx\\
& \leq L(M)(\|\tilde{u}-u\|_V + \|\tilde{a}-a\|_\calX)
\end{align*}
with $M= \left(C_{W\rightarrow W^{1,\frac{2p}{p-2}}} + C_{\calX\rightarrow L^\infty}\right)(\|u\|_V+\|\tilde{u}\|_V + \|c\|_\calX+\|\tilde{c}\|_\calX)$, subject to the constraint
\begin{flalign} \label{a-cond-2}
&\qquad\qquad\quad W= H^2(\Omega) \cap H_0^1(\Omega) \hookrightarrow L^{\frac{2p}{p-2}}(\Omega)\mbox{ i.e., } \quad p\geq\frac{2\bar{p}}{\bar{p}-2}.&
\end{flalign}

\item \textit{Well-definedness and boundedness of the parameter-to-state map:} \medbreak
A straightforward verification of 
{boundedness and coercivity} gives unique existence of the solution $u\in W(0,T)\subset \tilde{\calU}$ for $a\in\calD(F)\subset{\calX\hookrightarrow L^\infty(\Omega)}, \varphi\in L^2(0,T;H^{-1}(\Omega)), u_0\in L^2(\Omega)$. 

Similarly to the c-problem, the fact that the coercivity property of $f$ holds
\begin{align*}
\langle {-f(t,a,u)},u \rangle_{H^{-1},H_0^1}
&=\int_\Omega -\nabla\cdot(a\nabla u)u dx \geq \underline{a}\|u\|_{H_0^1(\Omega)}
\end{align*}
with the coefficient $\underline{a}$ being independent of $a$ shows uniform boundedness of $S$.

\item \textit{Continuous dependence on data of the solution to the linearized problem with zero initial data:} \medbreak

We employ the result in \cite[Section 7.1.3, Theorem 5]{Evans} with noting that the actual smoothness condition needed for the coefficient is that, $a$ is differentiable a.e on $\Omega$ and $a\in W^{1,\infty}(\Omega)$ rather than $a\in C^1(\Omega)$. From the observation $a\in\calD(F)=W^{1,p}(\Omega), p>d$ is differentiable a.e and the fact that $W^{1,\infty}(\Omega)$ is dense in $W^{1,p}(\Omega$), it enables us to imitate the contraction scenario and the dual argument as in the c-problem.


Taking $u,v$ solving (\ref{a-problem-1})-(\ref{a-problem-3}), we see
\begin{align*}
&\T: L^2(0,T;H^2(\Omega)\cap H_0^1(\Omega)) \rightarrow L^2(0,T;H^2(\Omega)\cap H_0^1(\Omega))\\
&\T=\left(\frac{d}{dt}-\nabla\cdot\Bigl(\hat{a}\nabla \Bigr)\right)^{-1}\nabla\cdot\Bigl((a-\hat{a})\nabla \Bigr)
\end{align*}
is a contraction
\begin{align} \label{a-contraction}
&\|\T( u -v)\|_{L^2(0,T; H^2\cap H_0^1)} \nonumber\\
&\leq \left\|\left(\frac{d}{dt}-\nabla\cdot\Bigl(\hat{a}\nabla \Bigr)\right)^{-1}\right\|_{L^2(0,T;L^2(\Omega))\rightarrow L^2(0,T;H^2\cap H_0^1)}
\hspace*{-2cm}\|\nabla\cdot\Bigl((a-\hat{a})\nabla(u-v)\Bigr)\|_{L^2(0,T;L^2(\Omega))}  \nonumber\\
&\leq C^{\hat{a}}\|\hat{a}-a\|_\calX \|u-v\|_{L^2(0,T;H^2\cap H_0^1)}\nonumber\\ 
&\leq C\epsilon \|u-v\|_{L^2(0,T;H^2\cap H_0^1)},
\end{align}
where $C\epsilon<1$ if we assume $\hat{a}=a^0\in W^{1,\infty}(\Omega)$ and $\rho$ is sufficiently small. If the index $p$ is large enough, smallness of $\rho$ can be omitted (discussed at the end of \ref{R3}).
Therefore, given $\phi\in L^2(0,T;L^2(\Omega))$, the adjoint state $p\in L^2(0,T; H^2\cap H_0^1)$ uniquely exists.

We also have the estimate
\begin{align*}
 \|p\|_{L^2(0,T;W)}
&\leq C \|\nabla\cdot\Bigl((a-\hat{a})\nabla p\Bigr)\|_{L^2(0,T;L^2(\Omega))}  + \|\phi\|_{L^2(0,T;L^2(\Omega))}) \\
&\leq C (2\rho\|p\|_{L^2(0,T;H^2\cap H_0^1)} + \|\phi\|_{L^2(0,T;L^2(\Omega))} )\\
&\leq C\|\phi\|_{L^2(0,T;V^*)},
\end{align*}
which proves continuous dependence of $p$ on $\phi\in L^2(0,T;V^*)$, consequently, continuous dependence of the solution $z\in L^2(0,T;V)$ on the data $b\in L^2(0,T;W^*)$ in \eqref{linear-eq-1}-\eqref{linear-eq-2}. Here smallness of $\rho$ is assumed.\\
If $p\geq 4$, smallness of $\rho$ is not required. To verify this, we test the adjoint equation by $-\Delta p$
\begin{align} \label{a-rho-1}
\int_\Omega-\dot{p}\Delta p + a(\Delta p)^2dx
&=\int_\Omega(-\nabla a\nabla p-\phi)\Delta p dx \nonumber\\
\frac{1}{2}\frac{d}{dt}\|\nabla p\|^2_{L^2(\Omega)}+\underline{a}\|\Delta p\|^2_{L^2(\Omega)}
&\leq \frac{\underline{a}}{2}\|\Delta p\|^2_{L^2(\Omega)} +\frac{1}{\underline{a}}\|\phi\|^2_{L^2(\Omega)}+\frac{1}{\underline{a}}\|\nabla a\nabla p\|^2_{L^2(\Omega)} \qquad\nonumber\\
\frac{1}{2}\frac{d}{dt}\|\nabla p\|^2_{L^2(\Omega)}+\frac{\underline{a}}{2}\|\Delta p\|^2_{L^2(\Omega)}
&\leq \frac{1}{\underline{a}}\|\phi\|^2_{L^2(\Omega)} + \frac{1}{\underline{a}}\|\nabla a\nabla p\|^2_{L^2(\Omega)},
\end{align}
where the last term on the right hand side can be estimated as in (\ref{c-coercive}) of the c-problem with $(\nabla a)^2$ in place of $c$, $\nabla p$ in place of $u$ and the assumption $\calX\hookrightarrow W^{1,4}(\Omega)$
\begin{align} \label{a-rho-2}
&\frac{1}{\underline{a}}\|\nabla a\nabla p\|^2_{L^2(\Omega)}\nonumber\\
&\leq \frac{C_{H_0^1\rightarrow L^6}}{\underline{a}}\|\nabla a\|^2_{L^4(\Omega)}\left( \frac{1}{16\epsilon\epsilon_1}\|\nabla p\|^2_{L^2(\Omega)} + \left(\frac{\epsilon_1}{4\epsilon} + \epsilon\right)\|\nabla p\|_{H_0^1(\Omega)}^2 \right) \nonumber\\
&\leq \frac{2C_{H_0^1\rightarrow L^6}}{\underline{a}}(\|a^0\|^2_{\calX}+\rho^2)\left( \frac{1}{16\epsilon\epsilon_1}\|\nabla p\|^2_{L^2(\Omega)} + \left(\frac{\epsilon_1}{4\epsilon} + \epsilon\right)\|\Delta p\|_{L^2(\Omega)}^2 \right).
\end{align}
Choosing $\epsilon_1<\epsilon$, and $\epsilon_1,\epsilon$ sufficiently small such that we can move the term involving $\|\Delta p\|_{L^2(\Omega)}^2$ from the right hand side to the left hand side of \eqref{a-rho-1}.
Note that, this choice of $\epsilon_1,\epsilon$ is just subject to $a^0$ and $\rho$.\\
Proceeding similarly to the c-problem, meaning applying Gronwall-Bellman inequality then taking the integral on $[0,T]$, we obtain
\begin{align} \label{a-rho-3}
\|p\|_{L^\infty(0,T;H^1(\Omega))}+\|\Delta p\|_{L^2(0,T;L^2(\Omega))}\leq C \|\phi\|^2_{L^2(0,T;L^2(\Omega))} 
\end{align}
with a constant $C$ depending only on $a^0, \rho$.

Observing the similarity in the form of the adjoint problem and the original problem, invoking the uniform bound \eqref{a-rho-3} w.r.t parameter $a$ and the fact $\overline{W^{1,\infty}(\Omega)}=W^{1,p}(\Omega)$ one can eliminate the need of smallness of $\rho$.

\item \textit{All-at-once tangential cone condition:} \medbreak
According to \eqref{qhat}, \eqref{stmnpq_aprob} with $s=0$, $t=2$, $m=n=2$, we require
\[
\frac{p}{p-1}\leq q\leq\hat{q}\geq2 \mbox{ and }1-\frac{d}{2}\succeq-\frac{d(p-1)}{p}+\frac{d}{q}\mbox{ and }-\frac{d}{2}\geq-d+\frac{d}{p}-1\,.\]
\end{enumerate}
\begin{corollary}
Assume $u_0\in L^2(\Omega), \varphi\in L^2(0,T;H^{-1}(\Omega))$, and
\begin{align} \label{a-space}
\begin{split}
\calX=W^{1,p}(\Omega),  \qquad Y=L^q(\Omega),\qquad &V=L^2(\Omega),  \qquad  W=H^2(\Omega)\cap H_0^1(\Omega) \\
& p\geq 2,\quad q\in\left[\underline{q},2\right],\quad d< p,
\end{split}
\end{align}
where 
$\underline{q}=\max\left\{\frac{p}{p-1}\,, \ \min\limits_{q\in[1,\infty]} \left\{1-\frac{d}{2}\succeq-\frac{d(p-1)}{p}+\frac{d}{q}\wedge-\frac{d}{2}\geq-d+\frac{d}{p}-1\right\}\right\}$.

Then $F$ defined by $F(a)=u$ solving \eqref{a-problem-1}-\eqref{a-problem-3} satisfies the tangential cone condition (\ref{TCC-Re}) with a uniformly bounded operator $F'(a)$ defined by (\ref{F'}).
\end{corollary}

\begin{remark}
This yields the possibility of a full Hilbert space setting $p=q=2$ of $X$ and $Y$ in case $d=1$, 
see also \cite{Hanke} and, for the static case, \cite{HNS95}.
\end{remark}

\subsection{An inverse source problem with a quadratic first order nonlinearity}

By the transformation $\utransf:=e^u$, the initial-value problem (\ref{log-problem-1})--(\ref{log-problem-3}) can be converted into an inverse potential problem as considered in Section \ref{c-problem}
\begin{alignat}{3} 
&\dot{\utransf}-\Delta \utransf +\theta\utransf=0 \qquad &&(t,x) \in (0,T)\times \Omega& \label{transf-log-problem-1}\\
& \utransf_{|\partial\Omega}=1 &&t \in (0,T)&	\label{transf-log-problem-2}\\
& \utransf(0)=\utransf_0 && x \in \Omega &	\label{transf-log-problem-3}
\end{alignat}
with $\utransf_0=e^{u_0}$.
Thus, in principle it is covered by the analysis from the previous section, as long as additionally positivity of $\utransf$ can be established. 
So the purpose of this section is to investigate whether we can allow for different function spaces $X,Y$ by directly considering (\ref{log-problem-1})--(\ref{log-problem-3}) instead of (\ref{transf-log-problem-1})--(\ref{transf-log-problem-3}).

We show that $f$ verifies the hypothesis proposed for the tangential cone condition in the reduced setting
on the function spaces
\begin{align} 
 \calX=L^p(\Omega),  \quad\qquad Y=L^q(\Omega),  \quad\qquad  V=W=H^2(\Omega)\cap H_0^1(\Omega).
\end{align}

\begin{enumerate}[label=(R\arabic*)]
\item \textit{Local Lipschitz continuity of $f$:}
\begin{align*}
\|-|\nabla\tilde{u}|^2 &+|\nabla u|^2 - \tilde{\theta}+\theta\|_{W^*}
= \sup_{\|w\|_W \leq 1} \int_\Omega \left(\nabla(u-\tilde{u})\cdot\nabla(u+\tilde{u})- \tilde{\theta} + \theta \right)w dx\\
&\leq C_{W\rightarrow L^{\bar{p}}}\left(\|(\nabla(u-\tilde{u})\cdot\nabla(u+\tilde{u})\|_{L^{\frac{\bar{p}}{\bar{p}-1}}} + \|\theta-\tilde{\theta}\|_{L^{\frac{\bar{p}}{\bar{p}-1}}} \right)\\
&\leq C_{W\rightarrow L^{\bar{p}}}\left(\|\nabla(u-\tilde{u})\|_{L^{\frac{2\bar{p}}{\bar{p}-1}}}\|\nabla(u+\tilde{u})\|_{L^{\frac{2\bar{p}}{\bar{p}-1}}} + \|\theta-\tilde{\theta}\|_{L^{\frac{\bar{p}}{\bar{p}-1}}} \right)\\
&\leq C_{W\rightarrow L^{\bar{p}}}\left( C^2_{V\rightarrow W^{1,{\frac{2\bar{p}}{\bar{p}-1}}}}\|u-\tilde{u}\|_V\|u+\tilde{u}\|_V + C_{\calX\rightarrow L^{\frac{\bar{p}}{\bar{p}-1}}} \|\theta-\tilde{\theta} \|_\calX \right).
\end{align*}
We can chose $L(M)=C_{W\rightarrow L^{\bar{p}}}\left(C^2_{V\rightarrow W^{1,{\frac{2\bar{p}}{\bar{p}-1}}}} (\|u\|_V+\|\tilde{u}\|_V ) + C_{\calX\rightarrow L^{\frac{\bar{p}}{\bar{p}-1}}}\right) +1$, under the conditions
\begin{equation} \label{log-cond-1}
\begin{split}
&V=H^2(\Omega)\cap H_0^1(\Omega) \hookrightarrow W^{1,{\frac{2\bar{p}}{\bar{p}-1}}}(\Omega)\mbox{ i.e., } \quad 1-\frac{d}{2}\geq -\frac{d(\bar{p}-1)}{2\bar{p}} \\
&\calX=L^p(\Omega) \hookrightarrow L^{\frac{\bar{p}}{\bar{p}-1}}(\Omega)\mbox{ i.e., }\quad p \geq \frac{\bar{p}}{\bar{p}-1}. 
\end{split}
\end{equation}

\item \textit{Well-definedness and boundedness of parameter-to-state map:}\medbreak

We argue unique existence of the solution to (\ref{log-problem-1})--(\ref{log-problem-3}) via the transformed problem 
(\ref{transf-log-problem-1})--(\ref{transf-log-problem-3}) for $\utransf=e^u$.

\quad To begin, by a similar argument to (\ref{c-contraction}) with the elliptic operator $A=-\Delta+\theta, \theta\in L^p(\Omega)$ in place of the parabolic operator, we show that the corresponding elliptic problem admits a unique solution in $H^2(\Omega)\cap H_0^1(\Omega)$ if the index $p$ satisfies (\ref{c-cond-3}). Employing next the semigroup theory in \cite[Section 7.4.3, Theorem 5]{Evans} or \cite[Chapter 7, Corollary 2.6]{Pazy} with assuming that $\utransf_0\in D(A)=H^2(\Omega)\cap H_0^1(\Omega)$ implies unique existence of a solution $\utransf\in C^1(0,T;H^2(\Omega))$ to \eqref{transf-log-problem-1}--\eqref{transf-log-problem-3}.
  
\quad Let $\utransf,\hat{\utransf}$ respectively solve \eqref{transf-log-problem-1}--\eqref{transf-log-problem-3} associated with the coefficients $\theta\in \calX, \hat{\theta}\in L^\infty(\Omega)$ with the same boundary and initial data, then $v=\utransf-\hat{\utransf}$ solves 
\begin{flalign*} 
&\qquad\qquad\qquad\dot{v}(t)-\Delta v(t) + \hat{\theta}v(t) = (\hat{\theta}-\theta)\utransf(t) \qquad t \in (0,T) &\\
&\qquad\qquad\qquad v(0)=0. &	
\end{flalign*}
Owing to the regularity from \cite[Section 7.1.3, Theorem 5]{Evans} and estimating similarly to (\ref{c-contraction}), we obtain
\begin{align} \label{log-positive-estimate-1}
\|\utransf-\hat{\utransf}\|_{L^\infty(0,T;H^2(\Omega))} 
&\leq C^{\hat{\theta}}\|(\hat{\theta}-\theta)\utransf\|_{H^1(0,T;L^2(\Omega))} \nonumber\\
&\leq C\|\hat{\theta}-\theta\|_\calX\|\utransf\|_{H^1(0,T;H^2(\Omega))}
\end{align}
with positive $\hat{\utransf}$ since $\hat{\theta}\in L^\infty(\Omega)$ and the constant $C$ depending only on $\theta^0, \rho$. Here we assume $\hat{\theta}=\theta^0\in L^\infty(\Omega)$ and $\rho$ is sufficiently small such that the right hand side is sufficiently small.
Then $\utransf\in L^\infty(0,T;H^2(\Omega))\subseteq L^\infty((0,T)\times\Omega)$ is close to $\hat{\utransf}$ and therefore positive as well. This assertion is valid if $0< \utransf_0=e^{u_0}\in H^2(\Omega)\cap H_0^1(\Omega)$, $0<\utransf|_{\delta\Omega}$, which is chosen as $\utransf|_{\delta\Omega}=1$ in this case (such that $\log(\utransf|_{\delta\Omega})=0$) and 
\begin{equation} \label{log-cond-2}
\begin{split}
&H^2(\Omega)\hookrightarrow L^{\frac{2p}{p-2}}(\Omega)\mbox{ i.e.,} \quad p\geq\frac{2\bar{p}}{\bar{p}-2} \quad\qquad\qquad\qquad\qquad\\
&V= H^2(\Omega)\hookrightarrow L^\infty(\Omega)\mbox{ i.e.,}\quad d\leq 3. 
\end{split}
\end{equation}

This leads to unique existence of the solution $u:=\log(\utransf)$ to the problem (\ref{log-problem-1})--(\ref{log-problem-3}), moreover  $0<\underline{c}\leq \utransf \in C^1(0,T;H^2(\Omega))$ allows $u=\log(\utransf) \in  C^1(0,T;H^2(\Omega)\cap H_0^1(\Omega))$.\\
If $d=1, p\geq2$, no assumption on smallness of $\rho$ is required since
\begin{align} \label{log-positive-estimate-2}
\|\utransf-\hat{\utransf}\|_{L^\infty(0,T;H^1(\Omega))}
&\leq C^\theta\|(\hat{\theta}-\theta)\hat{\utransf}\|_{L^2(0,T;L^2(\Omega))}\leq C\|\hat{\theta}-\theta\|_\calX\|\hat{\utransf}\|_{L^2(0,T;H^2(\Omega))} 
\end{align}
due to the estimates \eqref{c-rho-1}--\eqref{c-regularity} in Section \ref{c-problem}. Here the constant $C$ depends only on $\theta^0, \rho$ as claimed in (\ref{c-regularity}). This and the fact $\overline{L^\infty(\Omega)}=L^p(\Omega)$ allow us to chose $\hat{\theta}\in L^\infty(\Omega)$ being sufficiently close to $\theta\in L^p(\Omega)$ to make the right hand side of (\ref{log-positive-estimate-2}) arbitrarily small without the need of smallness of $\rho$.\\
\quad We have observed that, with the same positive boundary and initial data, the solution $\utransf=\utransf(\theta)$ to (\ref{transf-log-problem-1})--(\ref{transf-log-problem-3}) is bounded away from zero for all $\theta\in\calB^\calX_\rho(\theta^0)$. Besides, $S:\theta\mapsto \utransf$ is a bounded operator as proven in \ref{R2} of Section \ref{c-problem}. Consequently, $u=\log(\utransf)$ with $\Delta u=-\frac{|\nabla \utransf|^2}{\utransf^2}+\frac{\Delta \utransf}{\utransf}$ is uniformly bounded in $L^2(0,T;H^2(\Omega)\cap H_0^1(\Omega))$ for all $\theta\in\calB^\calX_\rho(\theta^0)$,  thus $S:\theta\mapsto u$ is a bounded operator on $\calB^\calX_\rho(\theta^0)$.

Moreover, we can derive a uniform bound for $\utransf$ in $H^1(0,T;H^2(\Omega))$ with respect to $\theta$. From
\begin{align*}
(\dot{\utransf}-\dot{\hat{\utransf}}) - \Delta(\utransf-\hat{\utransf}) + (\theta-\hat{\theta})(\utransf-\hat{\utransf})=-\hat{\theta}(\utransf-\hat{\utransf})-(\theta-\hat{\theta})\hat{\utransf},
\end{align*}
by taking the time derivative of both sides then test them with $-\Delta(\dot{\utransf}-\dot{\hat{\utransf}})$ we have
\begin{align*}
\frac{1}{2}\frac{d}{dt}\|\nabla(\dot{\utransf}-\dot{\hat{\utransf}})\|^2_{L^2(\Omega)}&+\|\Delta(\dot{\utransf}-\dot{\hat{\utransf}})\|^2_{L^2(\Omega)}\\
&\leq C_{H^2\hookrightarrow L^\infty}\|\theta-\hat{\theta}\|_{L^2(\Omega)}\|\Delta(\dot{\utransf}-\dot{\hat{\utransf}})\|^2_{L^2(\Omega)}\\
&\qquad+\|\hat{\theta}\|_{L^\infty(\Omega)}\|\dot{\utransf}-\dot{\hat{\utransf}}\|_{L^2(\Omega)}\|\Delta(\dot{\utransf}-\dot{\hat{\utransf}})\|_{L^2(\Omega)}\\
&\qquad +C_{H^2\hookrightarrow L^\infty}\|\theta-\hat{\theta}\|_{L^2(\Omega)}\|\Delta\dot{\hat{\utransf}}\|_{L^2(\Omega)}\|\Delta(\dot{\utransf}-\dot{\hat{\utransf}})\|_{L^2(\Omega))}\\
\frac{1}{2}\frac{d}{dt}\|\nabla(\dot{\utransf}-\dot{\hat{\utransf}})\|^2_{L^2(\Omega)}&+(1-\rho C_{H^2\hookrightarrow L^\infty}-\epsilon)\|\Delta(\dot{\utransf}-\dot{\hat{\utransf}})\|^2_{L^2(\Omega)}\\
&\leq \frac{1}{2\epsilon}\left(\|\hat{\theta}\|^2_{L^\infty(\Omega)}\|\dot{\utransf}-\dot{\hat{\utransf}}\|^2_{L^2(\Omega)}+C^2_{H^2\hookrightarrow L^\infty}\rho^2\|\Delta\dot{\hat{\utransf}}\|_{L^2(\Omega)}^2 \right),
\end{align*}
where $\|\Delta\dot{\hat{\utransf}}\|_{L^2(\Omega)}$ is attained by estimating with the same technique for \eqref{transf-log-problem-1}--\eqref{transf-log-problem-3} with the coefficient $\hat{\theta}\in L^\infty(\Omega)$. Since $\epsilon$ is arbitrarily small, if $\rho$ is sufficiently small and the following condition holds 
\begin{align}\label{log-cond-3}
\calX= L^p(\Omega)\hookrightarrow L^2(\Omega)\mbox{ i.e.,}\quad p\geq 2,
\end{align}
applying Gronwall's inequality then integrating on $[0,T]$ yields
\begin{equation}\label{log-positive-estimate-3}
\|\utransf-\hat{\utransf}\|_{H^1(0,T;H^2(\Omega))}\leq C\|\hat{\theta}-\theta\|_\calX\|\hat{\utransf}\|_{H^1(0,T;H^2(\Omega))} 
\end{equation}
for fixed $\hat{\utransf}=S(\hat{\theta})=S(\theta^0)$. So, $S(\calB^\calX_\rho(\theta^0))$ is bounded in $H^1(0,T;H^2(\Omega))$ and its diameter can be controlled by $\rho$. In case $d=1$, smallness of $\rho$ can be omitted if one uses the estimate \eqref{log-positive-estimate-2}.

\item \textit{Continuity of the inverse of the linearized model:}\medbreak

Now we consider the linearized problem
\begin{flalign} 
&\qquad\qquad\qquad\dot{z}(t)-\Delta z(t) + 2\nabla u(t)\cdot \nabla z(t) = r(t) \qquad t \in (0,T) & \label{log-linearized-1}\\
&\qquad\qquad\qquad z(0)=0,  & \label{log-linearized-2}	
\end{flalign}
whose adjoint problem after transforming $t=T-\tau$ is
\begin{flalign} 
&\qquad\qquad\qquad\dot{p}(t)-\Delta p(t) - 2\nabla\cdot \left(\nabla u(t) p(t)\right) = \phi(t) \qquad t \in (0,T) & \label{log-adjoint-1}\\
&\qquad\qquad\qquad p(0)=0 . \label{log-adjoint-2}& 	
\end{flalign}
Since $u\in C^1(0,T;H^2(\Omega)\cap H_0^1(\Omega))$ as proven in \ref{R2}, this equation with the coefficients $m:=-2\nabla u\in C^1(0,T;H^1(\Omega)), n:=-2\Delta u\in C^1(0,T;L^2(\Omega))$ is feasible to attain the estimate \ref{R3} by the contraction argument.

\quad Indeed, let us take $p$ solving (\ref{log-adjoint-1})--(\ref{log-adjoint-2}), 
then
\begin{flalign*}
\qquad \dot{p} &-\Delta p+\hat{m}\cdot\nabla p+\hat{n}p =\phi+ (\hat{m}-m)\cdot\nabla p+(\hat{n}-n)p&\\
\qquad p &= \left(\frac{d}{dt}-\Delta + \hat{m}\cdot\nabla+\hat{n}\right)^{-1}\left[\phi+(\hat{m}-m)\cdot\nabla p+(\hat{n}-n)p\right]&\\
&=: \T p&
\end{flalign*}
with some $\hat{m}\in L^\infty((0,T)\times\Omega)$ and some $\hat{n}\in L^\infty((0,T)\times\Omega)$ approximating $m$ and $n$. Then for $d\leq 3$, $\T: L^2(0,T;H^2(\Omega)\cap H_0^1(\Omega)) \rightarrow L^2(0,T;H^2(\Omega)\cap H_0^1(\Omega))$
is a contraction
\begin{align} \label{log-contraction}
&\|\T(p - q) \|_{L^2(0,T; H^2\cap H_0^1)} \nonumber\\
&\leq \left\|\left(\frac{d}{dt}-\Delta + \hat{m}\cdot\nabla+\hat{n}\right)^{-1}\right\|_{L^2(0,T;L^2(\Omega))\rightarrow L^2(0,T;H^2\cap H_0^1)}\nonumber\\
&\qquad .\left(\|(\hat{m}-m)\cdot\nabla(p-q)\|_{L^2(0,T;L^2(\Omega))}+\|(\hat{n}-n)(p-q)\|_{L^2(0,T;L^2(\Omega))} \right)\nonumber\\
&\leq C^{\hat{\theta}} \Big(\|\hat{m}-m\|_{L^\infty(0,T;H^1(\Omega))}\|\nabla(p-q)\|_{L^2(0,T;H^1(\Omega))} \nonumber\\
&\qquad\qquad\qquad\qquad\qquad\qquad+\|\hat{n}-n\|_{L^\infty(0,T;L^2(\Omega))}\|p-q\|_{L^2(0,T;L^\infty(\Omega))} \Big)\nonumber\\
&\leq C\epsilon \|p-q\|_{L^2(0,T;H^2\cap H_0^1)}, 
\end{align}
where $H_0^1(\Omega)\hookrightarrow L^6(\Omega), H^2(\Omega)\cap H_0^1(\Omega)\hookrightarrow L^\infty(\Omega)$ for $d\leq 3$. Above, we apply from \cite[Section 7.1.3 , Theorem 5]{Evans} the continuity of $\left(\frac{d}{dt}-\Delta + \hat{m}\cdot\nabla+\hat{n}\right)^{-1}$ with noting that, although the theorem is stated for time-independent coefficients, the proof reveals it is still applicable for $\hat{m}=\hat{m}(t,x), \hat{n}=\hat{n}(t,x)$ being bounded in time and space.\\
The above constant $C^{\hat{\theta}}$, which depends on $\hat{m}\in \nabla\cdot S(\calB^\calX_\rho(\theta^0))\cap L^\infty(0,T;L^\infty(\Omega))$, $\hat{n}\in \Delta S(\calB^\calX_\rho(\theta^0))\cap L^\infty(0,T;L^\infty(\Omega))$ can be bounded by some constant $C$ depending only on $S(\theta^0)$ and the diameter of $S(\calB^\calX_\rho(\theta^0))$ similarly to Sections \ref{c-problem} and \ref{a-problem} if choosing $\hat{\theta}=\theta^0$. In order to make $C\epsilon$ less than one, we require $\|\hat{m}-m\|_{L^\infty(0,T;H^1(\Omega))}$ and $\|\hat{n}-n\|_{L^\infty(0,T;L^2(\Omega))}$ to be sufficiently small. Those conditions turn out to be uniform boundedness of $\|\hat{\utransf}-\utransf\|_{L^\infty(0,T;H^2(\Omega))}$ (or the diameter of $S(\calB^\calX_\rho(\theta^0))$, which can be seen as smallness of $\rho$ as in \eqref{log-positive-estimate-3} since $H^1(0,T)\hookrightarrow L^\infty(0,T)$. 
From that, existence of the dual state $p\in L^2(0,T; H^2(\Omega)\cap H_0^1(\Omega))$ for given $\phi\in L^2(0,T;L^2(\Omega))$ is shown.

\quad Then \ref{R3-dual} follows without adding further constraints on $p$
\begin{align*}
 \|p&\|_{L^2(0,T;H^2(\Omega))}\\
&\leq C (\|(\hat{m}-m)\cdot\nabla p\|_{L^2(0,T;L^2(\Omega))}+\|(\hat{n}-n)p\|_{L^2(0,T;L^2(\Omega))} + \|\phi\|_{L^2(0,T;L^2(\Omega))} ) \\
&\leq C\|\phi\|_{L^2(0,T;L^2(\Omega))}
\end{align*}
with constant $C$ depending only on some fixed $\hat{m}, \hat{n}$ and the assumption on smallness of $\rho$.
Here with the $L^2$-norm on the right hand side, the maximum $q$ is limited by $\bar{q}=2$.

Observing that the problem (\ref{log-adjoint-1})--(\ref{log-adjoint-2}) has the form of the a-problem written in (\ref{a-rho-1}), with $\underline{a}=1, \nabla a=-2\nabla u(t) \in L^6(\Omega)$ and the additional term in the last line of the right hand side, namely,
\begin{align}
\frac{1}{\underline{a}}\|np\|^2_{L^2(\Omega)}&=\|\Delta up\|^2_{L^2(\Omega)}\leq \|\Delta u\|^2_{L^2(\Omega)}\|p\|^2_{L^\infty(\Omega)}\nonumber\\
&\leq C^2_{H_0^1\rightarrow L^\infty}\|\Delta u\|^2_{L^2(\Omega)}\|\nabla p\|^2_{L^2(\Omega)}\label{log-rho-1}
\end{align}
if the dimension $d=1$.\\
The solution $u=S(\theta)$ also lies in some ball in $C^1(0,T;H^2(\Omega)\cap H^1_0(\Omega))$ for all $\theta\in\calB^\calX_\rho(\theta^0)$, as in \ref{R2} we have shown boundedness of the operator $S$.\\
It allows us to evaluate analogously to (\ref{a-rho-1})--(\ref{a-rho-2}) with taking into account the additional term (\ref{log-rho-1}) to eventually get
\begin{align*}
\|\Delta p\|_{L^2(0,T;L^2(\Omega))}\leq C \|\phi\|^2_{L^2(0,T;L^2(\Omega))} 
\end{align*}
with the constant $C$ depending only on $\theta^0, \rho$. Hence, if $d=1$, $\rho$ is not required to be small.

\item \textit{All-at-once tangential cone condition:} \medbreak
According to \eqref{qhat}, \eqref{condPsi}  with $s=t=2$, $m=n=2$, $\hat{\gamma}=0, \rho=2$ this follows if
\[
\begin{aligned}
&2-\frac{d}{2}\succeq 1-\frac{d}{q^*}+\frac{d}{R}\mbox{ and }\\
& 1\leq\frac{R}{q^*}\mbox{ and } q\leq\hat{q} \mbox{ and }
2-\frac{d}{2}\succeq \max\left\{-\frac{d}{\hat{q}}\,, \ 1-\frac{d}{R}\right\},
\end{aligned}
\]
where the latter conditions come from the requirements $V=H^2(\Omega)\cap H_0^1(\Omega)\hookrightarrow W^{1,R}(\Omega)$.
\end{enumerate}

\begin{corollary}
Assume $u_0\in V$ and
\begin{align} \label{log-space}
\begin{split}
\calD(F)=\calX=L^p(\Omega),  \qquad Y=L^q(\Omega),\qquad &V=W=H^2(\Omega)\cap H_0^1(\Omega) \\ 
&p\geq2,\quad q\in\left[\underline{q},2\right],\quad d\leq 3
\end{split}
\end{align}
with $\underline{q}=\min\limits_q \left\{ 2-\frac{d}{2}\succeq 1-d+\frac{d}{q}+\frac{d}{\check{p}} \wedge q\geq 1+\frac{1}{\check{p}-1} \right\}$.

Then $F$  defined by $F(\theta)=u$ solving \eqref{log-problem-1}-\eqref{log-problem-3} satisfies the tangential cone condition (\ref{TCC-Re}) with a uniformly bounded operator $F'(\theta)$ defined by (\ref{F'}).
\end{corollary}

\begin{remark}
To achieve a Hilbert space setting for $X$ and $Y$, one can choose $p=q=2$ if $d\leq 3$,
see also \cite{Nguyen}.
\end{remark}

\subsection{An inverse source problem with a cubic zero order nonlinearity}\label{subsec:cubic}

We investigate this problem in the function spaces
\begin{flalign*} 
\calX=L^p(\Omega),  \qquad\qquad Y=L^q(\Omega),  \qquad\qquad  V=W=H_0^1(\Omega).&
\end{flalign*}

In the following we examine the conditions required for deriving the tangential cone condition and boundedness of the derivative of the forward operator.

\begin{enumerate}[label=(R\arabic*)]
\item \textit{Local Lipschitz continuity of $f$:}
\begin{align*}
&\|\tilde{u}^3- u^3 + \tilde{\theta}-\theta\|_{W^*}
= \sup_{\|w\|_W \leq 1} \int_\Omega (\tilde{u}-u)(\tilde{u}^2+\tilde{u}u+u^2)w + (\tilde{\theta}-\theta)w dx\\
&\leq C_{W\rightarrow L^{\bar{p}}}\left( \|(\tilde{u}-u)(\tilde{u}^2+\tilde{u}u+u^2)\|_{L^\frac{\bar{p}}{\bar{p}-1}} + \|\tilde{\theta}-\theta\|_{L^\frac{\bar{p}}{\bar{p}-1}} \right)\\
&\leq C_{W\rightarrow L^{\bar{p}}}\left(2 \|\tilde{u}-u\|_{L^{\bar{p}}}(\|\tilde{u}\|^2_{L^\frac{2\bar{p}}{\bar{p}-2}}+\|u\|^2_{L^\frac{2\bar{p}}{\bar{p}-2}})  + \|\tilde{\theta}-\theta\|_{L^\frac{\bar{p}}{\bar{p}-1}} \right)\\
&\leq C_{W\rightarrow L^{\bar{p}}}\left( 2C_{V\rightarrow L^{\bar{p}}}C^2_{V\rightarrow L^\frac{2\bar{p}}{\bar{p}-2}} \|\tilde{u}-u\|_V (\|\tilde{u}\|^2_V+\|u\|^2_V)  + \|\tilde{\theta}-\theta\|_\calX C_{\calX\rightarrow L^{\frac{\bar{p}}{\bar{p}-1}}}  \right). \qquad\quad
\end{align*}
We chose $L(M)= C_{W\rightarrow L^{\bar{p}}}\left( 2C_{V\rightarrow L^{\bar{p}}}C^2_{V\rightarrow L^\frac{2\bar{p}}{\bar{p}-2}} (\|\tilde{u}\|^2_V+\|u\|^2_V)  +  C_{\calX\rightarrow L^{\frac{\bar{p}}{\bar{p}-1}}}  \right)+1$, subject to the conditions
\begin{equation} \label{u3-cond-1}
\begin{split}
& 
V=W= H_0^1(\Omega) \hookrightarrow L^{\bar{p}}(\Omega)\mbox{ i.e.,} \quad 1-\frac{d}{2}\succeq-\frac{d}{\bar{p}}
\\
& V= H_0^1(\Omega) \hookrightarrow L^{\frac{2\bar{p}}{\bar{p}-2}}(\Omega)\mbox{ i.e.,} \quad d \leq 4 \qquad\qquad\qquad\qquad\qquad\qquad\\
& \calX=L^p(\Omega) \hookrightarrow L^{\frac{\bar{p}}{\bar{p}-1}}(\Omega)\mbox{ i.e.,}\quad p \geq \frac{\bar{p}}{\bar{p}-1}.
\end{split}
\end{equation}

\item \textit{Well-definedness and boundedness of the parameter-to-state map:}\medbreak
Verifying the conditions \ref{S1}-\ref{S3} 
with the Gelfand triple $H_0^1(\Omega)\hookrightarrow L^2(\Omega)\hookrightarrow H^{-1}(\Omega)$  
shows that the problem (\ref{u3-problem-1})-(\ref{u3-problem-3}) admits a unique solution in the space $W(0,T)$. Subsequently, \cite[Theorem 8.16]{Roubicek} strengthens the solution to belong to $ L^\infty(0,T;V)$. To validate this regularity result, the following additional assumptions are made
\begin{align} \label{u3-cond-2}
& \calX=L^p(\Omega) \hookrightarrow L^2(\Omega)\mbox{ i.e., }\quad p \geq 2, &
\end{align}
the initial data $u_0\in V$ and the known source term $\varphi\in L^2(0,T;L^2(\Omega))$.

From \cite[Proposition 4.2, Section 6.1]{Nguyen}, we have
\begin{align*}
\|S(\theta)&\|_{L^\infty(0,T;V)}
\leq N\left(\|\theta+\varphi\|_{L^2(0,T;L^2(\Omega)}+\sqrt{\int_\Omega\frac{1}{2}|\nabla u_0|^2+\frac{1}{4}u_0^4 dx } \right)\\
&\leq N\left(\sqrt{T}(\|\theta_0\|_{L^2(\Omega)}+\rho)+\|\varphi\|_{L^2(0,T;L^2(\Omega)}+\sqrt{\int_\Omega\frac{1}{2}|\nabla u_0|^2+\frac{1}{4}u_0^4 dx } \right)
\end{align*}
for some $N$ depending only on $c^\theta_0=c_0=\frac{1}{2}$.
This thus implies uniform boundedness of $S$ on $\calB^\calX_\rho(\theta_0)$.

\item \textit{Continuous dependence on data of the solution to the linearized problem with zero initial data:}\medbreak

For this purpose, semi-coercivity of the linearized forward operator is obvious
\begin{align*}
\dupair{-f'_u(t,\theta,v),v}_{V^*,V} &= \int_\Omega (-\Delta v + 3u^2v)v dx		\\
&\geq \|\nabla v\|^2_{L^2(\Omega)}=\|v\|^2_V.
\end{align*}

\item \textit{All-at-once tangential cone condition:} \medbreak 
According to \eqref{qhat}, \eqref{condPhi}, with $s=t=1$, $m=n=2$, $\gamma=\kappa=1, r=\hat{q}=\bar{p}$ this follows if
\[
2\leq \frac{\bar{p}}{q^*}\mbox{ and }1-\frac{d}{2}\succeq-\frac{d}{q^*}+\frac{2d}{\bar{p}}\mbox{ and }
 q\leq\bar{p} \mbox{ and }1-\frac{d}{2}\succeq -\frac{d}{\bar{p}},
\]
where the latter condition comes from the requirement $V=H_0^1(\Omega)\hookrightarrow L^{\bar{p}}(\Omega).$
\end{enumerate}

\begin{corollary}
Assume $u_0\in H_0^1(\Omega), \varphi\in L^2(0,T;L^2(\Omega))$, and
\begin{align} \label{u3-space}
\begin{split}
\calD(F)=\calX=L^p(\Omega),  \qquad\qquad &Y=L^q(\Omega),  \qquad\qquad  V=W=H_0^1(\Omega)\qquad\\
& p\geq 2,\quad q\in \left[\underline{q}, \bar{q}\right],\quad d\leq 4,
\end{split}
\end{align}
where $\underline{q}=\min\limits_q \left\{1-\frac{d}{2}\succeq-d+\frac{d}{q}+\frac{2d}{\bar{p}} \wedge q\geq 1+\frac{2}{\bar{p}-2} \right\}$ with
\begin{align}
d=1 \text{ and } \bar{q}=\infty,\qquad
d=2 \text{ and } \bar{q}<\infty,\qquad
d\geq 3 \text{ and } \bar{q}=\frac{2d}{d-2}. 
\end{align}

Then $F$ defined by $F(\theta)=u$ solving \eqref{u3-problem-1}-\eqref{u3-problem-3} satisfies the tangential cone condition (\ref{TCC-Re}) with a uniformly bounded operator $F'(\theta)$ defined by (\ref{F'}).
\end{corollary}

\begin{remark}
Here $X$ and $Y$ can be chosen as Hilbert spaces with $p=q=2$ and $d\leq 3$.
\end{remark}

\section*{Appendix A}
\numberwithin{equation}{section}
\renewcommand{\thesection}{A}
\renewcommand{\theequation}{\thesection.\arabic{equation}}

\subsection*{Notation}
\begin{enumerate}[label=$\bullet$]
\item For $a,b\in\R$, the notation $a\succeq b$ means: $a\geq b$ with strict inequality if $b=0$. 
\item For normed spaces $A,B$, the notation $A\hookrightarrow B$ means: $A$ is continuously embedded in $B$.
\item For a normed space $A$, an element $a\in  A$ and $\rho>0$, we denote by $\calB_\rho^A(a)$ the closed ball of radius $\rho$ around $a$ in $A$.
\item For vectors $\vec{a},\vec{b}\in \R^n$, $\vec{a}\cdot\vec{b}$ denotes the Euclidean inner product. Likewise, $\nabla\cdot \vec{v}$ denotes the divergence of the vector field $\vec{v}$. 
\item $C$ denotes a generic constant that may take different values whenever it appears.
\item For $p\in[1,\infty]$, we denote by $p^*=\frac{p}{p-1}$ the dual index.
\end{enumerate}

\subsection*{proof of Proposition \ref{prop:cprob_aao}}
On \eqref{Xcprob}
for some $p\in[1,\infty]$, we can estimate by applying H\"older's inequality, once with exponent $p$ and once with exponent $\frac{q}{p^*}$ (where $p^*=\frac{p}{p-1}$ is the dual index)
\[
\begin{aligned}
&\|(B\hat{c})(t)\hat{v}\|_{W^*}
=\sup_{w\in W\,, \ \|w\|_W\leq1}\int_\Omega \hat{c}\, \hat{v}\, w\, dx
\leq \|\hat{c}\|_{L^p} \|\hat{v}\|_{L^q} \sup_{w\in W\,, \ \|w\|_W\leq1} \|w\|_{L^{\frac{p^* q}{q-p^*}}},
\end{aligned}
\]
where we need to impose $q\geq p^*$ and in case of equality formally set $\frac{p^* q}{q-p^*}=\infty$. In order to guarantee continuity of the embedding $W\hookrightarrow L^{\frac{p^* q}{q-p^*}}(\Omega)$ as needed here, we therefore, together with \eqref{qhat}, require the conditions
\begin{equation}\label{stmnpq_cprob}
s-\frac{d}{m}\succeq -\frac{d}{\hat{q}}\mbox{ and }\hat{q}\geq q\geq p^* \mbox{ and } t-\frac{d}{n}\succeq  -\frac{d(q-p^*)}{p^*q}\,.
\end{equation}

\subsection*{proof of Proposition \ref{prop:aprob_aao}}

With $X$ as in \eqref{Xaprob}, in order to guarantee the required boundedness of the embeddings 
\[
\begin{aligned}
&\mathcal{X}\hookrightarrow L^r(\Omega)\,, \quad
W\hookrightarrow W^{1,\frac{p^* q}{q-p^*}}(\Omega)\,, \
W\hookrightarrow W^{2,\frac{r^* q}{q-r^*}}(\Omega)\,, \\
&\mbox{for some $r\in[1,\infty]$ such that $r^*\leq q$}
\end{aligned}
\]
we impose, additionally to \eqref{qhat}, the conditions
\[
\begin{aligned}
&(a) \ \hat{q}\geq q\geq \max\{p^*,r^*\} \mbox{ and } 
&& (b) \ t-1-\frac{d}{n}\succeq-\frac{d(q-p^*)}{p^*q} \mbox{ and } \\
& (c) \ t-2-\frac{d}{n}\succeq-\frac{d(q-r^*)}{r^*q} \mbox{ and } 
&& (d) \ 1-\frac{d}{p}\succeq-\frac{d}{r}
\end{aligned}
\]
for some $r\in[1,\infty]$. To eliminate $r$, observe that the requirement $(c)$, i.e., 
$t-2-\frac{d}{n}\succeq-\frac{d}{r^*}+\frac{d}{q}$ gets weakest when $r^*$ is chosen minimal, which, subject to requirement (d) is 
\begin{equation}\label{rstar}
r\begin{cases}=\infty\mbox{ if }p>d\\<\infty\mbox{ if }p=d\\=\frac{dp}{d-p}\mbox{ if }p<d\end{cases} \quad\mbox{, i.e., } \quad
r^*\begin{cases}=1\mbox{ if }p>d\\>1\mbox{ if }p=d\\=\frac{dp}{dp-d+p}\mbox{ if }p<d\end{cases}\,.
\end{equation}
Inserting this into (c) and taking into account \eqref{qhat}, we end up with the following requirements on $s,t,p,q,m,n$ (using the fact that $q\geq p^*$ implies $q\geq\frac{dp}{dp-d+p}$):
\begin{equation}\label{stmnpq_aprob}
\begin{aligned}
&s-\frac{d}{m}\succeq -\frac{d}{\hat{q}}\mbox{ and }
\hat{q}\geq q \geq p^* \mbox{ and }\\ 
&t-1-\frac{d}{n}\succeq-\frac{d(q-p^*)}{p^*q} \mbox{ and }
t-2-\frac{d}{n}\begin{cases}\succeq-d+\frac{d}{q}\mbox{ if } p>d\\>-d+\frac{d}{q}\mbox{ and }q>1\mbox{ if } p=d\\
\succeq-\frac{dp-d+p}{p}+\frac{d}{q}\mbox{ if } p<d\,.\end{cases}
\end{aligned}
\end{equation}

\subsection*{proof of Proposition \ref{prop:nlinvsource_aao}}

Here we have 
\[
\begin{aligned}
&\Bigl(f(\theta,u)-f(\ttheta,\tu)-f_u'(\theta,u)(u-\tu)-f_\theta'(\theta,u)(\theta-\ttheta)\Bigr)(t)\\
&=
\int_0^1\Bigl(\Phi'((u(t)+\sigma(\tu(t)-u(t))-\Phi'(u(t))\Bigr)\,d\sigma\,(\tu(t)-u(t))\\
&\quad +\int_0^1\Bigl(\Psi'(\nabla u(t)+\sigma(\nabla \tu(t)-\nabla u(t))-\Psi'(\nabla u(t)) \Bigr)\,d\sigma\,\nabla(\tu(t)-u(t))\,.
\end{aligned}
\]
This shows that the only condition which has to be taken into account when choosing the space $\mathcal{X}$ is that $B(t)\in \cL(\mathcal{X},W^*)$.
Again we assume $\mathcal{C}(t)$ to be the embedding operator $V\hookrightarrow Y$.

As opposed to Section \ref{subsec:bilin}, where we could do the estimates pointwise in time, we will now also have to to use H\"older estimates with respect to time.
To this end, we dispose over the following continuous embeddings
\[
\begin{aligned}
&\mathcal{U}\hookrightarrow L^2(0,T;W^{s,m}(\Omega))\\
&\mathcal{U}\hookrightarrow L^\infty(0,T;H^{\ts}(\Omega)) \mbox{ provided }W^{s-\ts,m}(\Omega)\hookrightarrow W^{t+\ts,n}(\Omega)\,,
\end{aligned}
\]
where the first holds just by definition of $\mathcal{U}$ and the second follows from \cite[Lemma 7.3]{Roubicek}\footnote{$L^2(0,T;\widetilde{W})\cap H^1(0,T;\widetilde{W}^*)\hookrightarrow L^\infty(0,T;L^2(\Omega))$} with $\widetilde{W}=W^{t+\ts,n}(\Omega)$, using the fact that
\[
\begin{aligned}
&u\in L^2(0,T;W^{s,m}(\Omega))\cap H^1(0,T;(W^{t,n}(\Omega))^*)\\
& \Leftrightarrow \ D^{\ts} u\in L^2(0,T;W^{s-\ts,m}(\Omega))\cap H^1(0,T;(W^{t+\ts,n}(\Omega))^*),
\end{aligned}
\]
where $D^\ts v=\sum_{|\alpha|\leq\ts} D^\alpha v $.

We first consider the case of an affinely linear (or just vanishing) function $\Psi$, which still comprises, e.g., models with linear drift and diffusion, so that $C_{\Psi''}$ can be set to zero. We can then estimate 
\[
\begin{aligned}
&\|f(\theta,u)-f(\ttheta,\tu)-f_u'(\theta,u)(u-\tu)-f_\theta'(\theta,u)(\theta-\ttheta)\|_{L^2(0,T;W^*)}\\
&\leq C_{\Phi''} \left(\int_0^T\, \Bigl(\sup_{w\in W\,, \ \|w\|_W\leq1}\int_\Omega 
(1+|u(t)|^{\gamma}+|\tu(t)|^{\gamma})\,|\tu(t)-u(t)|^{1+\kappa} \, w\, dx\Bigr)^2\, dt\right)^{1/2}\,, 
\end{aligned}
\]
where, using H\"older's inequality three times ($P=q$, $P=\frac{r}{q^*(\gamma+\kappa)}$, $P=\frac{\gamma+\kappa}{\gamma}$) and continuity of the embedding $H^{\ts}(\Omega)\hookrightarrow L^r(\Omega)$ provided $\ts-\frac{d}{2}\succeq-\frac{d}{r}$
\begin{equation}\label{est_Phi}
\begin{aligned}
&\left(\int_0^T\, \Bigl(\sup_{w\in W\,, \ \|w\|_W\leq1}\int_\Omega |u(t)|^{\gamma}\,|\tu(t)-u(t)|^{1+\kappa} \, w\, dx\Bigr)^2\, dt\right)^{1/2}\\
&\leq \|\tu-u\|_{L^2(0,T;L^q(\Omega))} 
\sup_{w\in W\,, \ \|w\|_W\leq1}\Bigl\||u|^{\gamma} |\tu-u|^\kappa w\Bigr\|_{L^\infty(0,T;L^{q^*}(\Omega))}\\
&\leq \|\tu-u\|_{\mathcal{Y}} \Bigl\|\Bigl(|u|^{\gamma} |\tu-u|^\kappa\Bigr)^{\frac{1}{\gamma+\kappa}}\Bigr\|_{L^\infty(0,T;L^r(\Omega))}^{\gamma+\kappa}
\sup_{w\in W\,, \ \|w\|_W\leq1}\| w\|_{L^{\frac{rq^*}{r-q^*(\gamma+\kappa)}}(\Omega)}\\
&\leq \|\tu-u\|_{\mathcal{Y}} \|u\|_{L^\infty(0,T;L^r(\Omega))}^{\gamma} \|\tu-u\|_{L^\infty(0,T;L^r(\Omega))}^\kappa
\sup_{w\in W\,, \ \|w\|_W\leq1}\| w\|_{L^{\frac{rq^*}{r-q^*(\gamma+\kappa)}}(\Omega)}\\
&\leq (C_{H^{\ts}\to L^r}^\Omega)^{\gamma+\kappa} \|u\|_{L^\infty(0,T;H^{\ts}(\Omega))}^{\gamma} \|\tu-u\|_{L^\infty(0,T;H^{\ts}(\Omega))}^\kappa\\
&\qquad \|\tu-u\|_{\mathcal{Y}} \sup_{w\in W\,, \ \|w\|_W\leq1}\| w\|_{L^{\frac{rq^*}{r-q^*(\gamma+\kappa)}}(\Omega)}
\end{aligned}
\end{equation}
(and likewise for the term containing $|\tu(t)|^{\gamma}$)
for some $r\in[1,\infty]$ with $\frac{r}{q^*}\geq\gamma+\kappa$.
In order to get finiteness of the $L^\infty(0,T;H^{\ts}(\Omega))$ norms appearing here by means of \cite[Lemma 7.3]{Roubicek}, we assume the embedding $W^{s-{\ts},m}(\Omega)\hookrightarrow W^{t+{\ts},n}(\Omega)$ to be continuous, which leads to the condition
\[
s-\ts-\frac{d}{m}\succeq t+\ts-\frac{d}{n}\mbox{ and }s-\ts\geq t+\ts\,.
\]
Moreover, in order to guarantee continuity of the embedding $W\hookrightarrow L^{\frac{rq^*}{r-q^*(\gamma+\kappa)}}(\Omega)$ and for the above H\"older estimate to make sense we impose
\[
\gamma+\kappa\leq \frac{r}{q^*}\mbox{ and }t-\frac{d}{n}\succeq-\frac{d(r-q^*(\gamma+\kappa))}{rq^*}
\]
for some $r\in[1,\infty]$.
Summarizing, we have the following conditions
\begin{equation}\label{summcond}
\begin{aligned}
&\ts-\frac{d}{2}\succeq-\frac{d}{r}\mbox{ and }
s-\ts-\frac{d}{m}\succeq t+\ts-\frac{d}{n}\mbox{ and }s-\ts\geq t+\ts\mbox{ and }\\
&\gamma+\kappa\leq \frac{r}{q^*}\mbox{ and }t-\frac{d}{n}\succeq-\frac{d(r-q^*(\gamma+\kappa))}{rq^*}=-\frac{d}{q^*}+\frac{d(\gamma+\kappa)}{r},
\end{aligned}
\end{equation}
which imply 
\[
s\geq \frac{d}{m}+d-\frac{d}{q^*}+d\frac{\gamma+\kappa-2}{r}\,.
\]
This lower bound on $s$ gets weakest for maximal $r$, if $\gamma+\kappa>2$ and for minimal $r$ if $\gamma+\kappa<2$. We therefore make the following case distinction.\\
If $\gamma+\kappa>2$ or $\gamma+\kappa=2$ and $q=1$ we set $r=\infty$, which leads to $\ts>\frac{d}{2}$, hence, according to \eqref{summcond}, we can choose
\begin{equation}\label{casegammkappage2}
\begin{aligned}
&\mbox{case  $\gamma+\kappa>2$ or ($\gamma+\kappa=2$ and $q=1$):}\\
&t>\frac{d}{n}-\frac{d}{q^*}\,, \quad q\leq \hat{q}\,, \\
&s>\max\left\{t+d+\max\left\{0,\frac{d}{m}-\frac{d}{n}\right\},\frac{d}{m}-\frac{d}{\hat{q}}\right\}.
\end{aligned}
\end{equation}
If $\gamma+\kappa<2$ or $\gamma+\kappa=2$ and $q>1$ we set $r=\max\{1,q^*(\gamma+\kappa)\} <\infty$, $\ts:=\max\{0,\frac{d}{2}-\frac{d}{r}\}$ and, according to \eqref{summcond}, can therefore choose
\begin{equation}\label{casegammkappale2}
\begin{aligned}
&\mbox{case  $\gamma+\kappa<2$ or ($\gamma+\kappa=2$ and $q>1$):}\\
&t>\frac{d}{n} +\min\left\{0,-\frac{d}{q^*}+d(\gamma+\kappa)\right\}\,, \quad q\leq \hat{q}\,,\\ 
&s>\max\left\{t+\max\left\{0,d-\frac{2d}{\max\{1,q^*(\gamma+\kappa)\}}\right\},\frac{d}{m}-\frac{d}{\hat{q}}\right\}\,.
\end{aligned}
\end{equation}

\medskip

Now we consider the situation of nonvanishing gradient nonlinearites $C_{\Psi''}>0$ where we additionally need to estimate terms of the form 
\[
\left(\int_0^T\, \Bigl(\sup_{w\in W\,, \ \|w\|_W\leq1}\int_\Omega |\nabla u(t)|^{\hat{\gamma}}\,|\nabla\tu(t)-\nabla u(t)|^{1+\hat{\kappa}} \, w\, dx\Bigr)^2\, dt\right)^{1/2}\,,
\]
which, in  order to end up with an estimate in terms of $\|\tu-u\|_{L^2(0,T;L^q(\Omega))}$ requires us to move the gradient by means of integration by parts. Assuming for simplicity that $\hat{\kappa}=1$ we get 
\[
\begin{aligned}
&\left(\int_0^T\, \Bigl(\sup_{w\in W\,, \ \|w\|_W\leq1}\int_\Omega |\nabla u(t)|^{\hat{\gamma}}\,|\nabla\tu(t)-\nabla u(t)|^2 \, w\, dx\Bigr)^2\, dt\right)^{1/2}\\
&=\left(\int_0^T\, \Bigl(\sup_{w\in W\,, \ \|w\|_W\leq1}\int_\Omega (\tu(t)-u(t))\, g^w(t) \, dx\Bigr)^2\, dt\right)^{1/2}\\
&\leq \|\tu-u\|_{L^2(0,T;L^q(\Omega))} 
\sup_{w\in W\,, \ \|w\|_W\leq1}\|g^w\|_{L^\infty(0,T;L^{q^*}(\Omega))},
\end{aligned}
\]
where 
\[
\begin{aligned}
g^w(t)&=\nabla\cdot\Bigl(|\nabla u(t)|^{\hat{\gamma}}\,\nabla(\tu(t)-u(t)) \, w\Bigr)\\
&=\hat{\gamma}|\nabla u(t)|^{\hat{\gamma}-2}(\nabla^2 u(t) \nabla u(t))\cdot\nabla(\tu(t)-u(t)) \, w\\
&\quad+|\nabla u(t)|^{\hat{\gamma}}\,\Delta(\tu(t)-u(t)) \, w
+|\nabla u(t)|^{\hat{\gamma}}\,\nabla(\tu(t)-u(t)) \cdot\nabla w\\
&\quad=:g_1(t)+g_2(t)+g_3(t),
\end{aligned}
\]
where $\nabla^2$ denotes the Hessian.
For the last term we proceed analogously to above (basically replacing $u$ by $\nabla u$ and $w$ by $\nabla w$) to obtain 
\begin{equation}\label{est_Psi1}
\begin{aligned}
&\|g_3\|_{L^\infty(0,T;L^{q^*}(\Omega))}=\|\, |\nabla u(t)|^{\hat{\gamma}}\,\nabla(\tu(t)-u(t)) \cdot\nabla w\|_{L^\infty(0,T;L^{q^*}(\Omega))}\\
&\leq \|\nabla u\|_{L^\infty(0,T;L^R(\Omega))}^{\hat{\gamma}} \|\nabla(\tu-u)\|_{L^\infty(0,T;L^R(\Omega))}
\sup_{w\in W\,, \ \|w\|_W\leq1}\|\nabla w\|_{L^{\frac{Rq^*}{R-q^*(\hat{\gamma}+1)}}(\Omega)}
\end{aligned}
\end{equation}
and use \cite[Lemma 7.3]{Roubicek} with $\nabla u\in L^2(0,T;W^{s-1,m}(\Omega))\cap H^1(0,T;(W^{t+1,n}(\Omega))^*)$, which under the conditions 
\begin{equation}\label{condts}
\begin{aligned}
&t-\frac{d}{n}\succeq 1-\frac{d(R-q^*(\hat{\gamma}+1))}{Rq^*}\,,\\
&s-1-\ts-\frac{d}{m}\succeq t+1+\ts-\frac{d}{n}\,, \ s-1-\ts\geq t+1+\ts\,, \ \ts-\frac{d}{2}\succeq-\frac{d}{R} 
\end{aligned}
\end{equation}
yields $\nabla u\in L^\infty(0,T;H^\ts(\Omega))\subseteq L^\infty(0,T;L^R(\Omega))$ and $W\hookrightarrow W^{1,\frac{Rq^*}{R-q^*(\hat{\gamma}+1)}}(\Omega)$.

The other two terms can be bounded by 
\[
|g_1(t)+g_2(t)|\leq 
\Bigl(\hat{\gamma} |\nabla^2 u(t)| \, |\nabla u(t)|^{\hat{\gamma}-1}\, |\nabla(\tu(t)-u(t))|
+|\nabla^2 (\tu(t)-u(t))| \, |\nabla u(t)|^{\hat{\gamma}}\Bigr) \,|w|
\]
(note that here $|\cdot|$ denotes the Frobenius norm of a matrix)
so that it suffices to find an estimate on expressions of the form
\[
\||\nabla^2 z| \, |\nabla v|^{\hat{\gamma}-1} \, |\nabla y| \,|w|\|_{L^\infty(0,T;L^2(\Omega))}
\]
for $z,v,y\in \mathcal{U}$, $w\in W$.
To this end, we will again employ \cite[Lemma 7.3]{Roubicek}, making use of the fact that for any $\varrho ,R\in[1,\infty)$, due to H\"older's inequality with $P=\frac{\varrho }{2}$ and with $P=\frac{R(\varrho -2)}{2\varrho \hat{\gamma}}$, the estimate 
\begin{equation}\label{est_Psi2}
\begin{aligned}
&\|\,|\nabla^2 z| \, |\nabla v|^{\hat{\gamma}-1} |\nabla y| \, |w|\|_{L^2(\Omega)}\\
&\leq \|\,|\nabla^2 z| \,\|_{L^\varrho (\Omega)}
\|\Bigl(|\nabla v|^{\hat{\gamma}-1} |\nabla y|\Bigr)^{\frac{1}{\hat{\gamma}}}\|_{L^R(\Omega)}^{\hat{\gamma}}
\|w\|_{L^{\frac{2R\varrho }{R(\varrho -2)-2\varrho \hat{\gamma}}}(\Omega)}\\
&\leq C_{H^{\tts}\to L^\varrho }^\Omega (C_{H^{\ttts}\to L^\varrho }^\Omega )^{\hat{\gamma}}
\|\, |\nabla^2 z|\|_{H^\tts(\Omega)} \|\Bigl(|\nabla v|^{\hat{\gamma}-1} |\nabla y|\Bigr)^{\frac{1}{\hat{\gamma}}}\|_{H^\ttts(\Omega)}^{\hat{\gamma}} 
\|w\|_{L^{\frac{2R\varrho }{R(\varrho -2)-2\varrho \hat{\gamma}}}(\Omega)}
\end{aligned}
\end{equation}
holds.
To make sense of these H\"older estimates and to guarantee continuity of the embedding $W\hookrightarrow L^{\frac{2R\varrho }{R(\varrho -2)-2\varrho \hat{\gamma}}}(\Omega)$ we impose
\begin{equation}\label{rR}
\begin{aligned}
\varrho \geq 2\mbox{ and }R\geq\frac{2\varrho \hat{\gamma}}{\varrho -2}\mbox{ and } 
t-\frac{d}{n}\succeq -\frac{d(R(\varrho -2)-2\varrho \hat{\gamma})}{2R\varrho }=-\frac{d}{2}+\frac{d}{\varrho }+\frac{d\hat{\gamma}}{R}
\end{aligned}
\end{equation}
Taking into account the fact that here $\nabla^2 z$ contains second and $\Bigl(|\nabla v|^{\hat{\gamma}-1} |\nabla y|\Bigr)^{\frac{1}{\hat{\gamma}}}$ first derivatives of elements of $\mathcal{U}$, we therefore aim at continuity of the embeddings
\[
\begin{aligned}
&L^2(0,T;W^{s-2,m}(\Omega))\cap H^1(0,T;W^{t+2,n}(\Omega)) \hookrightarrow L^\infty(0,T;H^\tts(\Omega)) \hookrightarrow L^\infty(0,T;L^\varrho (\Omega))\\
&L^2(0,T;W^{s-1,m}(\Omega))\cap H^1(0,T;W^{t+1,n}(\Omega)) \hookrightarrow L^\infty(0,T;H^\ttts(\Omega)) \hookrightarrow L^\infty(0,T;L^R(\Omega))\,,
\end{aligned}
\]
which can be achieved by means of \cite[Lemma 7.3]{Roubicek} under the conditions
\begin{equation}\label{ttsttts}
\begin{aligned}
&s-2-\tts-\frac{d}{m}\succeq t+2+\tts-\frac{d}{n}\mbox{ and }s-2-\tts\geq t+2+\tts\mbox{ and }\tts-\frac{d}{2}\succeq-\frac{d}{\varrho }\\
&s-1-\ttts-\frac{d}{m}\succeq t+1+\ttts-\frac{d}{n}\mbox{ and }s-1-\ttts\geq t+1+\ttts\mbox{ and }\ttts-\frac{d}{2}\succeq-\frac{d}{R}\,.
\end{aligned}
\end{equation}

For instance, we may set $\varrho =2$, $R=\infty$ to obtain, inserting into \eqref{condts}, \eqref{rR}, \eqref{ttsttts}, that $\tts\geq0$, $\ttts>\frac{d}{2}$ hence 
\begin{equation}\label{condgrad}
\begin{aligned}
&t>\frac{d}{n}\,, \ t-\frac{d}{n}\succeq 1-\frac{d}{q^*}\,, \ s\succeq t+2+\max\{2,d\}+\frac{d}{m}-\frac{d}{n}\,, \ s\geq t+2+\max\{2,d\} \,, \\ 
&s -\frac{d}{m}\succeq -\frac{d}{\hat{q}}\,,  \ q\leq \hat{q}\,.
\end{aligned}
\end{equation}

\medskip

In order to avoid the use of too high values of $s$ and $t$, we can alternatively skip the use of \cite[Lemma 7.3]{Roubicek} and instead set 
\begin{equation}\label{U1}
\mathcal{U}=\{u \in L^\infty(0,T;L^r(\Omega))\cap L^2(0,T;V): \dot{u} \in L^2(0,T;W^*)\}
\end{equation}
in case $C_{\Psi''}=0$, or 
\begin{equation}\label{U2}
\mathcal{U}=\{u \in L^\infty(0,T;L^r(\Omega)\cap W^{1,R}(\Omega)\cap W^{2,\varrho}(\Omega))\cap L^2(0,T;V): \dot{u} \in L^2(0,T;W^*)\}
\end{equation}
otherwise.
This can also be embedded in a Hilbert space setting by replacing $L^\infty(0,T)$ with $H^\sigma(0,T)$ for some $\sigma>\frac12$.
Going back to estimate \eqref{est_Phi} in case $C_{\Phi''}=0$ we end up with the conditions
\begin{equation}\label{condPhi}
\gamma+\kappa\leq \frac{r}{q^*}\mbox{ and }t-\frac{d}{n}\succeq-\frac{d}{q^*}+\frac{d(\gamma+\kappa)}{r}\,,
\end{equation}
cf. \eqref{summcond},
and in case $C_{\Psi''}>0$, considering estimates \eqref{est_Psi1}, \eqref{est_Psi2} otherwise, we require
\begin{equation}\label{condPsi}
\begin{aligned}
&t-\frac{d}{n}\succeq \max\left\{1-\frac{d}{q^*}+\frac{d(\hat{\gamma}+1)}{R}, -\frac{d}{2}+\frac{d}{\varrho }+\frac{d\hat{\gamma}}{R}\right\}\mbox{ and }\\
& \hat{\gamma}+1\leq\frac{R}{q^*}\mbox{ and }\varrho \geq 2\mbox{ and }\hat{\gamma}\leq\frac{R(\varrho -2)}{2\varrho}\,,
\end{aligned}
\end{equation}
cf. \eqref{condts} \eqref{rR}, and in both cases we additionally need to impose \eqref{qhat}.

\section*{Acknowledgment}
BK is supported by the Austrian Science Fund (FWF) with project P30054
(Solving Inverse Problems without Forward Operators). 
OS is supported by the Austrian Science Fund (FWF) with project
F6807-N36 (Tomography with Uncertainties) and with project I3661-N27
(Novel Error Measures and Source Conditions of Regularization Methods
for Inverse Problems).

BK and OS acknowledge the support of BIRS for a stay at the Banff
center, Canada, where the paper has been finished.

\bibliographystyle{siam}
\nocite{*}
\bibliography{Reference}

\end{document}